\documentclass[preprint,sort&compress,final,12pt]{elsarticle}

\usepackage{amsmath}
\usepackage{amssymb}
\usepackage{graphicx}
\usepackage{latexsym}
\usepackage{epstopdf}
\usepackage{natbib}
\usepackage{color}
\usepackage{hyperref}
\usepackage{geometry}
\usepackage{bm}
\usepackage{amsthm}
\usepackage{amscd}
\usepackage{extarrows}
\usepackage{setspace}
\usepackage{mathrsfs}
\usepackage{subcaption}

\newtheorem{theorem}{Theorem}[section]

\newtheorem{lemma}[theorem]{Lemma}

\numberwithin{equation}{section}

\def\dsum{\displaystyle\sum}

\def\dfrac#1#2{\frac{\displaystyle {#1}}{\displaystyle {#2}}}

\newcommand{\Id}{\mathrm{Id}}

\journal{\empty}

\begin{document}
	
	\begin{frontmatter}
		
		\title{Bounded and unbounded solutions of asymmetric oscillations in the completely degenerate resonant case}
		
		\author[au1]{Yaqi Liang}
		
		\ead[au1]{202331130015@mail.bnu.edu.cn}
		
		\author[au2]{Xiong Li\footnote{Partially supported by the NSFC (12371158). Corresponding author.}}
		
		\ead[au2]{xli@bnu.edu.cn}
		
		\address[au1] {School of Science, Hebei University of Technology, Tianjin 300130, People's Republic of China}
		
		\address[au2]{Laboratory of Mathematics and Complex Systems (Ministry of Education), School of Mathematical Sciences, Beijing Normal University, Beijing 100875, People's Republic of China}

		\begin{abstract}
			This paper is concerned with the asymmetric oscillator $x^{\prime\prime}+a x^{+}-bx^{-}=f(t)$, where $x^{+}=\max\{x,0\}$, $x^{-}=\max\{-x,0\}$, \(f\) is real analytic \(2\pi\)-periodic, \(a\) and \(b\) are two distinct positive constants satisfying the resonance condition $\frac{1}{2}\left(\frac{1}{\sqrt{a}}+\frac{1}{\sqrt{b}}\right)=\frac{n}{m}.$ The discriminant function is 
			\[
			\Phi_{f}(\theta)
			=
			\int_{0}^{2\pi}
			C\left(\omega_{0}\theta+nt\right)f(nt)\,dt,
			\]
			where \(C(t)\) is the solution of $x^{\prime\prime}+a x^{+}-b x^{-}=0$ with $C(0)=1, C'(0)=0$. Li and Li \cite{unbounded} proved the equation admits unbounded solutions whenever $\Phi_f(\theta)\not\equiv0$ and the zero set of \(\Phi_f\) is nonempty. However, this result does not cover $\Phi_f(\theta)\equiv0$. Moreover, Li and Li \cite{unbounded} pointed out that, for \(a=4\), \(b=1\), and \(f(t)=\cos kt\), \(\Phi_f(\theta)\equiv0\) whenever \(4\nmid k\), and left the corresponding boundedness problem open. Thus we will investigate the completely degenerate case and resolve the above problem, and finally obtain a complete classification: the equation $x^{\prime\prime}+4 x^{+}- x^{-}=\cos kt$ admits unbounded solutions if \(k\in\{1,2\}\) or \(4\mid k\), whereas all its solutions are bounded if \(k\geq3\) and \(4\nmid k\).
		\end{abstract}
		
		\begin{keyword}
			Asymmetric oscillators; Poincar\'{e} maps; Unbounded solutions; Boundedness
		\end{keyword}
		
	\end{frontmatter}

	\section{Introduction}
	In the 1960s, Littlewood \cite{Littlewood1968} raised the problem of determining
	whether every solution of the periodically forced nonlinear oscillator
	\begin{equation*}
		x''+g(x)=p(t)
	\end{equation*}
	is bounded. Here, $p$ is $2\pi$-periodic and $g\colon\mathbb R\to\mathbb R$
	is continuous, with
	\[
	\operatorname{sign}(x)g(x)\longrightarrow +\infty
	\quad\text{as } \quad |x|\longrightarrow+\infty.
	\]
	Morris \cite{Morris1976} obtained the first fundamental boundedness result in
	the superlinear setting. Using a twist theorem, he proved that all solutions
	are bounded when $g(x)=2x^3$ and $p$ is piecewise continuous. Further
	developments for superlinear equations can be found in
	\cite{Dieckerhoff, You1992,Yuan1998} and the references therein.
	
	The semilinear case is considerably more delicate. If
	\[
	0<a\leq \frac{g(x)}{x}\leq b<+\infty,
	\qquad x\neq0,
	\]
	then the resonance may produce unbounded motions; the elementary equation
	\[
	x''+n^2x=\cos nt,
	\qquad n\in\mathbb N,
	\]
	has no bounded solution. At the same time, nonlinear effects may restore
	boundedness even in resonant regimes. For boundedness and periodic-solution
	results concerning semilinear Duffing equations, we refer to
	\cite{Ding1993,JiaoPiaoWang2012,LazerLeach1969}.
	
	In this paper we study the asymmetric oscillator
	\begin{equation}\label{equ}
		x''+ax^+-bx^-=f(t),
	\end{equation}
	where $x^+=\max\{x,0\}$, $x^-=\max\{-x,0\}$, $a$ and $b$ are distinct
	positive constants, and $f$ is a real-analytic $2\pi$-periodic function.
	Equations of this type arise in models of suspension bridges \cite{A.P} and
	belong to the class of equations with jumping nonlinearities introduced in the
	work of Fu\v{c}ik \cite{Fucik} and Dancer \cite{Dancer}; see also
	\cite{Gallouet,Habet,Lazer,Zhang}.
	
	Set
	\[
	\omega_0:=\frac12\left(\frac1{\sqrt a}+\frac1{\sqrt b}\right),
	\]
	and the least positive period of the corresponding unforced oscillator is
	$2\pi\omega_0$. The boundedness problem for \eqref{equ} was first addressed
	by Ortega \cite{Ortega1998}. In particular, for
	\[
	x''+ax^+-bx^-=1+\varepsilon h(t),
	\]
	with $h$ smooth and $2\pi$-periodic, he proved that all solutions are bounded
	provided $|\varepsilon|$ is sufficiently small. This phenomenon sharply
	contrasts with linear resonance: if $a=b=n^2$, then every solution of
	$x''+n^2x=1+\varepsilon\cos nt$ is unbounded whenever $\varepsilon\neq0$.
	
	The arithmetic character of $\omega_0$ separates the problem into
	nonresonant and resonant cases. When $\omega_0\notin\mathbb Q$, Ortega
	\cite{Ortega2001} proved that a nonzero mean of $f$ implies boundedness of all
	solutions. For zero-mean forcing, boundedness was subsequently established
	under a Bruno arithmetic condition in \cite{LiLi2021}, by means of the
	normal-form method of \cite{Russmann2002}. An approximate-normal-form
	criterion for zero-mean forcing was obtained in \cite{LiangLi2025}.
	
	We focus on the resonant case. Write $\omega_0=n/m$, where $m,n\in\mathbb N$
	are relatively prime. Let $C$ be the solution of
	\[
	x''+ax^+-bx^-=0,
	\qquad C(0)=1,\quad C'(0)=0.
	\]
	The function $C(t)$ is periodic with least positive period $2\pi\omega_0$. On the interval $[-\pi\omega_0,\pi\omega_0]$, it is given by
	\begin{equation*}
		C(t)=\left\{\begin{array}{l}
			\cos(\sqrt{a}\,|t|), \quad 0\leq |t|\leq \frac{\pi}{2\sqrt{a}}, \\[0.25cm]
			-\sqrt{\frac{a}{b}}\sin\left(\sqrt{b}\left(|t|-\frac{\pi}{2\sqrt{a}}\right)\right),\quad \frac{\pi}{2\sqrt{a}}\le|t|\leq \omega_0\pi.\end{array}
		\right.
	\end{equation*}
	Let $S=C^\prime$. Then $S(t)\in C(\mathbb{R})$, and the following properties hold:
	
	\vspace{1em}
	(1) $C(-t)=C(t),S(-t)=-S(t)$;
	
	\vspace{1em}
	(2) $C(t)$ and $S(t)$ are $2\pi\omega_0$-periodic functions;
	
	\vspace{1em}
	(3) $S(t)^2+aC^{+}(t)^{2}+bC^{-}(t)^{2}\equiv a$.
	
	\vspace{1em}
	The first resonant discriminant is
	\[
	\Phi_f(\theta)
	:=\int_0^{2\pi}
	C\left(\frac nm\theta+nt\right)f(nt)\,dt,
	\]
	and denote its zero set by
	\[
	\mathcal A(f):=\{\theta\in\mathbb R:\Phi_f(\theta)=0\}.
	\]
	Liu \cite{Liu} proved that all solutions of \eqref{equ} are bounded when
	$\mathcal A(f)=\varnothing$. Alonso and Ortega \cite{Alonso} showed that if
	every zero of $\Phi_f$ is nondegenerate, then all sufficiently large solutions
	are unbounded; their argument also implies the existence of an unbounded
	solution whenever $\Phi_f$ has at least one nondegenerate zero.
	
	The degenerate regime has proved substantially more difficult. Ortega
	\cite{Ortega1998} considered
	\[
	x''+4x^+-x^-=\lambda+\cos4t
	\]
	and left open the threshold cases $|\lambda|=1/45$, in which all zeros of the
	first resonant discriminant are degenerate. These cases were later resolved by Li and Li
	\cite{LiLi2022}, who proved the existence of unbounded solutions. Li and Li
	\cite{unbounded} subsequently showed, for analytic forcing, that an unbounded
	solution exists whenever $\Phi_f(\theta)\not\equiv0$ and
	$\mathcal A(f)\neq\varnothing$. These results do not cover the completely
	degenerate case $\Phi_f(\theta)\equiv0$. The following examples proposed by Li and Li \cite{unbounded} show that this situation can indeed happen. Consider the equation
	\begin{equation}\label{eq:model-intro}
		x''+4x^+-x^-=\cos kt,
	\end{equation}
	where \(k\) is a positive integer. Then the first resonant discriminant function is
	\[
	\Phi_f(\theta)
	=
	\begin{cases}
		0,
		& k=1,2,\\[6pt]
		\displaystyle
		-\frac{4}{(k^2-4)(k^2-1)}
		\bigl(1+(-1)^k\bigr)
		\left[
		\cos\left(\frac{k\pi}{4}+\frac{3k}{4}\theta\right)
		+
		\cos\left(\frac{3k\pi}{4}+\frac{3k}{4}\theta\right)
		\right],
		& k\geq 3.
	\end{cases}
	\]
	Thus, if \(k=4r\), \(r=1,2,3,\ldots\),
	\[
	\Phi_f(\theta)
	=
	\frac{4(-1)^{r+1}}
	{(4r^2-1)(16r^2-1)}
	\cos 3r\theta,
	\]
	and if \(4\nmid k\),
	\[
	\Phi_f(\theta)\equiv 0.
	\]
	
	Hence, by the result of Li and Li \cite{unbounded}, equation \eqref{eq:model-intro} admits an unbounded solution when $4\mid k$. However, for $ 4\nmid k$, it remained unknown whether all solutions of the equation \eqref{eq:model-intro} are bounded. This motivates us to investigate the completely degenerate case $\Phi_f(\theta)\equiv 0$. In this paper, we will first establish several sufficient conditions on \(f\) for the existence of  unbounded solutions when \(\Phi_f(\theta)\equiv0\), and then study the boundedness problem of equation \eqref{eq:model-intro}.
	
	Before stating the main theorem, we introduce some notations. Set
	\[
	\rho=\sqrt{\frac{2m}{an}},\quad\kappa=\rho m\left(\frac{n}{m}\right)^{\frac32},\quad\beta=\frac{n^2\rho^2}{2m},\quad\delta=\frac{\rho^3m}{8}\left(\frac{n}{m}\right)^{\frac52}.
	\]
	Define
	\begin{equation}\label{QQ}
		Q(\theta,t)
		=\int_{0}^{t}C(ns)f(\theta+ns)\,ds-\frac{1}{2\pi}\int_{0}^{2\pi}\left(\int_{0}^{\tau}C(ns)f(\theta+ns)\,ds\right)d\tau,
	\end{equation}
	\begin{equation*}
		G_0(\theta,t)=
		\beta C(nt)^2f(\theta+nt)^2-\frac{\kappa^2}{2}C(nt)f(\theta+nt)\partial_\theta Q(\theta,t),
	\end{equation*}
	\begin{equation*}
		G_1(\theta,t)
		=
		\delta C(nt)^3f(\theta+nt)^3
		-
		\frac{\kappa^3}{8}C(nt)f(\theta+nt)(\partial_\theta Q(\theta,t))^2,
	\end{equation*}
	\begin{equation*}
		G_2(\theta,t)=-\frac{\kappa^4}{16}C(nt)f(\theta+nt)(\partial_\theta Q(\theta,t))^3+\frac{\kappa\delta}{2}C(nt)^3f(\theta+nt)^3\partial_\theta Q(\theta,t),
	\end{equation*}
	and
	\begin{equation*}
		P_1(\theta,t)=G_0(\theta,t), \qquad [P_1](\theta)=\frac{1}{2\pi}\int_{0}^{2\pi}P_1(\theta,t)dt,
	\end{equation*}
	\begin{equation}\label{P_2}
		P_2(\theta,t)
		=
		G_1(\theta,t)+\frac{\kappa Q(\theta,t)}{2}\partial_\theta G_0(\theta,t), \qquad [P_2](\theta)=\frac{1}{2\pi}\int_{0}^{2\pi}P_2(\theta,t)dt,
	\end{equation}
	\begin{equation}\label{P_3}
		P_3(\theta,t)=
		G_2(\theta,t)
		+\frac{\kappa Q(\theta,t)}{2}\partial_\theta G_1(\theta,t)
		+\frac{\kappa^2Q(\theta,t)\partial_\theta Q(\theta,t)}{4}\partial_\theta G_0(\theta,t)
		+\frac{\kappa^2Q(\theta,t)^2}{8}\partial_\theta^2G_0(\theta,t),
	\end{equation}
	\begin{equation*}
		[P_3](\theta)=\frac{1}{2\pi}\int_{0}^{2\pi}P_3(\theta,t)dt.
	\end{equation*}
	
	Now we are in a position to introduce our main result.
	\begin{theorem}\label{th}
		Suppose that \(f\) is a real-analytic \(2\pi\)-periodic function and that \(\Phi_f(\theta)\equiv0\). Then the following conclusions hold:
		
		(1) if $[P_1]$ is nonconstant, then equation \eqref{equ} admits an unbounded solution;
		
		(2) if $[P_1]$ is constant and $[P_2]$ has a nondegenerate zero, then equation \eqref{equ} admits an unbounded solution;
		
		(3) if \([P_1]\) is constant, \([P_2](\theta)\equiv0\), and \([P_3]\) has a nondegenerate zero, then equation \eqref{equ} admits an unbounded solution.
	\end{theorem}
	
	Our approach is to continue the resonant normal form beyond the case where the first discriminant $\Phi_f$ vanishes. Building on the normal-form analysis of \cite{unbounded}, we carry the canonical normalization to three additional orders. After a sequence of canonical transformations, the Hamiltonian is reduced, in the high-energy region, to
	\[
	H(r,\theta,t)
	=[P_1](\theta)+[P_2](\theta)r^{-1/2}
	+[P_3](\theta)r^{-1}+O(r^{-3/2}).
	\]
	
	First, if the function \([P_1]\) defined above is nonconstant, then \eqref{equ} admits an unbounded solution. More precisely, when $[P_1]$ is nonconstant, the resulting system has a different asymptotic structure from that considered in~\cite{unbounded}: instead of constructing an invariant region, we establish a quantitative iteration lemma which keeps the angular variable in a suitably shrinking high-energy strip and forces the radial action to grow monotonically. This requires a delicate balance between the leading drift and the higher-order remainders.
	
	Second, if $[P_1]$ is a constant and $[P_2]$ has a nondegenerate zero, such a two-dimensional trapping argument is no longer effective, since the escaping orbit is confined to a one-dimensional invariant manifold. Indeed, the original Poincar\'{e} map possesses an \(n\)-resonant elliptic fixed point at infinity. After passing to the \(n\)-fold resonant return map, compactifying infinity via \(u=r^{-1/2}\), and applying a singular time rescaling, the resulting degenerate equilibrium is desingularized into a hyperbolic saddle. A Lyapunov--Perron contraction then selects the stable or unstable manifold producing the desired unbounded orbit. The same construction remains applicable after the further degeneration $[P_2](\theta)\equiv0$, provided that $[P_3]$ has a nondegenerate zero, with the desingularizing time scale adjusted accordingly.
	
	As an application of Theorem~\ref{th}, we consider equation \eqref{eq:model-intro} in the case $4\nmid k$, since it admits an unbounded solution when $4\mid k$ by \cite{unbounded}.
	\begin{theorem}\label{th2}
		For $k=1$ and $k=2$, equation \eqref{eq:model-intro} admits an unbounded solution. For $k\geq3$ with $4\nmid k$, all solutions of \eqref{eq:model-intro} are bounded.
	\end{theorem}
	
	For \(k=1\) and \(k=2\), these conclusions follow from parts (3) and (1) of Theorem \ref{th}, respectively. Indeed, a direct computation shows that \([P_1]\) is nonconstant when \(k=2\). When \(k=1\), \([P_1]\) is constant and \([P_2](\theta)\equiv0\); moreover, \([P_3]\) has a nondegenerate zero.
	
	For $4\nmid k$ and $k\geq 3$, the result is based on Lemma \ref{periodic}.
	In fact, by considering the successor map of equation \eqref{eq:model-intro}, Lemma \ref{periodic} shows that, in the high-energy region, the phase shifts accumulated during the motions on the two half-axes satisfy an exact arithmetic cancellation. As a result, a finite iterate of the successor map is the identity, so that the successor map in this region is periodic. Accordingly, all the nontrivial averaged coefficients in the successive normal forms vanish identically. This is not merely a finite-order degeneracy of the asymptotic expansion, but an exact degeneracy reflecting the periodic structure of the original piecewise-linear equation.
	
	Previous studies on bounded and unbounded solutions of asymmetric oscillators usually concentrated on the relationship between the global internal frequency $1/\omega_0$ and the forcing frequency of $f$. However, since an asymmetric oscillator has different restoring coefficients on the positive and negative half-axes, its motions on the two half-axes are governed by the local natural frequencies \(\sqrt{a}\) and \(\sqrt{b}\), respectively. To the best of our knowledge, this is the first result showing that resonance between the forcing frequency and the local natural frequency associated with either half-axis may produce an unbounded solution. In particular, for \(a=4\), \(b=1\), and \(f(t)=\cos kt\), the forcing frequency \(k=2\) resonates with the natural frequency on the positive half-axis, whereas \(k=1\) resonates with that on the negative half-axis; in both cases, the equation admits unbounded solutions. By contrast, for all the remaining forcing frequencies for which the first discriminant is completely degenerate, an exact phase cancellation ensures that all solutions remain bounded. Therefore, this paper reveals a half-axis resonance mechanism that cannot be detected by the first discriminant and thereby establishes a complete and refined frequency-selection law for this family of asymmetric oscillators.
	
	The paper is organized as follows. Section 2 derives the normal form and the
	coefficients $[P_1]$, $[P_2]$, and $[P_3]$. Section 3 treats the case in which
	$[P_1]$ is nonconstant. Section 4 develops the Lyapunov--Perron construction
	when $[P_1]$ is a constant and $[P_2]$ has a nondegenerate zero, while Section 5
	handles the further degenerate case in which $[P_1]$ is constant and $[P_2](\theta)\equiv0$. The final section applies the general theory to the equation \eqref{eq:model-intro} and establishes
	the bounded and unbounded behaviors described above.

	\section{The normal form}
	This section derives a higher-order normal form for equation \eqref{equ} under the completely degenerate resonance condition \(\Phi_f(\theta)\equiv0\). We first introduce action-angle variables, interchange the roles of the time and angle variables, and perform appropriate scaling and rotational transformations. We then successively eliminate the time-dependent parts of the first three higher-order terms by canonical averaging transformations. This procedure reduces the system in the high-energy region to the normal form \eqref{simple}, whose averaged coefficients \([P_1]\), \([P_2]\), and \([P_3]\) determine the subsequent analysis of unbounded solutions.
	
	\medskip
	\noindent\textbf{Step 1: Action-angle variable transformation.}
	Let $y=x^\prime$. Then equation \eqref{equ} is equivalent to the following planar system
	\begin{equation}\label{equ3}
		\left\{\begin{array}{l}
			x^\prime=y, \\ [0.25cm]
			y^\prime=-ax^{+}+bx^{-}+f(t). \end{array}
		\right.
	\end{equation}
	
	For $r>0$ and $\theta\in\mathbb{T}:=\mathbb{R}/2\pi\mathbb{Z}$, define the following polar-coordinate transformation $\Gamma :(r,\theta)\rightarrow (x,y)$:
	\[ 	\left\{\begin{array}{l}
		x=\rho r^{\frac{1}{2}}C\left(\frac{n}{m}\theta\right), \\ [0.25cm]
		y=\rho r^{\frac{1}{2}}S\left(\frac{n}{m}\theta\right), \end{array}
	\right.	\] where $\rho=\sqrt{\frac{2m}{an}}$ is a positive constant, and $r$ and $\theta$ are the action and angle variables, respectively. A direct calculation shows that $\Gamma$ is a symplectic transformation.
	
	In Cartesian coordinates, the Hamiltonian associated with system \eqref{equ3} can be written as \[H(x,y,t)=\frac{1}{2}y^2+\frac{a}{2}(x^+)^2+\frac{b}{2}(x^-)^2-f(t)x,\]
	and in the new coordinates $(r,\theta)$ it becomes
	\begin{equation}\label{equ4}
		H(r,\theta,t)=\frac{m}{n}r-\rho r^{\frac{1}{2}}C\left(\frac{n}{m}\theta\right)f(t).
	\end{equation}
	Consequently, system \eqref{equ3} becomes
	\begin{equation*}
		\left\{\begin{array}{l}
			\theta^\prime=H_r=\frac{m}{n}-\frac{1}{2}\rho r^{-\frac{1}{2}}C\left(\frac{n}{m}\theta\right)f(t),\\ [0.25cm]
			r^\prime=-H_\theta=\frac{n}{m}\rho r^{\frac{1}{2}}S\left(\frac{n}{m}\theta\right)f(t). \end{array}
		\right.
	\end{equation*}
	This is a semilinear system.
	
	We now interchange the roles of the time and angle variables; that is, we regard $\theta$ as the new time variable and $t$ as the new angle variable. The new action variable is then $H$, and the new Hamiltonian is $r=r(t,H,\theta)$. The function $r=r(t,H,\theta)$ can be solved from \eqref{equ4}. Indeed, \eqref{equ4} gives
	\begin{equation}\label{taylor}
		r^{\frac{1}{2}}=\frac{\rho Cf+\sqrt{\rho^{2}C^{2}f^{2}+4mH/n}}{2m/n}.
	\end{equation}
	Using the Taylor expansion, we may further expand \eqref{taylor} for sufficiently large $H$ to obtain
	\[
	\begin{aligned}
		r(H,t,\theta)
		&=\frac{n}{m}H+\rho\left(\frac{n}{m}\right)^{\frac32}
		C\left(\frac{n}{m}\theta\right)f(t)H^{\frac12}\\
		&+
		\frac{\rho^2}{2}\left(\frac{n}{m}\right)^2
		C^2\left(\frac{n}{m}\theta\right)f^2(t)\\
		&+
		\frac{\rho^3}{8}\left(\frac{n}{m}\right)^{\frac52}
		C^3\left(\frac{n}{m}\theta\right)f^3(t)H^{-\frac12}\\
		&-
		\frac{\rho^5}{128}\left(\frac{n}{m}\right)^{\frac72}
		C^5\left(\frac{n}{m}\theta\right)f^5(t)H^{-\frac32}+O(H^{-\frac52}).
	\end{aligned}
	\]
	The new Hamiltonian system is now
	\[
	\frac{dH}{d\theta}=-\frac{\partial r}{\partial t}(H,t,\theta),\qquad \frac{dt}{d\theta}=\frac{\partial r}{\partial H}(H,t,\theta),	\]
	which is still a semilinear system.
	
	\medskip
	\noindent\textbf{Step 2: Scaling transformation.} We next perform a scaling transformation. Introduce the new time variable $\tilde{\theta}=\frac{\theta}{m}$. Then
	\begin{equation*}	
		\left\{\begin{array}{l}
			\frac{dH}{d\tilde{\theta}}=-mr_t=-\frac{\partial}{\partial t}\tilde{r}(H,t,\tilde{\theta}),  \\ [0.25cm]
			\frac{dt}{d\tilde{\theta}}=mr_H=\frac{\partial}{\partial H}\tilde{r}(H,t,\tilde{\theta}), \end{array}
		\right.	
	\end{equation*}
	where
	\begin{equation*}
		\begin{aligned}
			\tilde{r}&=mr=nH+m\rho\left(\frac{n}{m}\right)^\frac{3}{2} CfH^\frac{1}{2}+\frac{n^2}{2m}\rho^2C^2f^2\\
			&+\frac{\rho^3m}{8}\left(\frac{n}{m}\right)^\frac{5}{2}C^3f^3H^{-\frac{1}{2}}-\frac{\rho^5m}{128}\left(\frac{n}{m}\right)^\frac{7}{2}C^5f^5H^{-\frac{3}{2}}+O\left(H^{-\frac{5}{2}}\right).
		\end{aligned}
	\end{equation*}
	For convenience, relabel $(\tilde{\theta},\tilde{r},H,t)$ as $(t,H,r,\theta)$. The new system can then be rewritten as
	\begin{equation*}
		\left\{\begin{array}{l}
			\frac{dr}{dt}=-\frac{\partial}{\partial\theta}H(r,\theta,t),  \\ [0.25cm]
			\frac{d\theta}{dt}=\frac{\partial}{\partial r}H(r,\theta,t),  \end{array}
		\right.
	\end{equation*}
	where
	\begin{equation}\label{equ6}
		\begin{aligned}
			H&=nr+m\rho\left(\frac{n}{m}\right)^\frac{3}{2} C(nt)f(\theta)r^\frac{1}{2}+\frac{n^2}{2m}\rho^2C^2(nt)f^2(\theta)\\
			&+\frac{\rho^3m}{8}\left(\frac{n}{m}\right)^\frac{5}{2}C^3(nt)f^3(\theta)r^{-\frac{1}{2}}-\frac{\rho^5m}{128}\left(\frac{n}{m}\right)^\frac{7}{2}C^5(nt)f^5(\theta)r^{-\frac{3}{2}}+O\left(r^{-\frac{5}{2}}\right).
		\end{aligned}
	\end{equation}
	
	\medskip
	\noindent\textbf{Step 3: Rotational transformation.} We now carry out canonical transformations to compute the normal form of the Hamiltonian system \eqref{equ6}. First, introduce a rotational transformation to eliminate the linear part of the Hamiltonian in \eqref{equ6}. Define the generating function of the rotational transformation by $S(r,t)=-nrt$. Then the rotational transformation $\Phi_1: (r_1,\theta_1,t)\rightarrow(r,\theta,t)$ can be written as
	\[ 	\left\{\begin{array}{l}
		r=r_1+\frac{\partial S}{\partial \theta}=r_1,  \\ [0.25cm]
		\theta=\theta_1-\frac{\partial S}{\partial r}=\theta_1+nt. \end{array}
	\right.	\]
	
	Under the transformation $\Phi_1$, the new Hamiltonian $H_1(r_1,\theta_1,t)$ is
	\begin{equation}\label{eq22}
		\begin{split}
			H_1(r_1,\theta_1,t)&=H+\frac{\partial S}{\partial t} \\
			&=m\left(\frac{n}{m}\right)^\frac{3}{2}\rho C(nt)f(\theta_1+nt)r_1^\frac{1}{2}+\frac{n^2}{2m}\rho^2C^2(nt)f^2(\theta_1+nt)\\
			&+\frac{\rho^3m}{8}\left(\frac{n}{m}\right)^\frac{5}{2}C^3(nt)f^3(\theta_1+nt)r_1^{-\frac{1}{2}}-\frac{\rho^5m}{128}\left(\frac{n}{m}\right)^\frac{7}{2}C^5(nt)f^5(\theta_1+nt)r_1^{-\frac{3}{2}}+O(r_1^{-\frac{5}{2}}).
		\end{split}
	\end{equation}
	Under $\Phi_1$, the semilinear system defined by the Hamiltonian in \eqref{equ6} is transformed into the sublinear system defined by the Hamiltonian in \eqref{eq22}. It is readily verified that the transformation preserves the periodicity and boundedness of solutions.
	
	\medskip
	\noindent\textbf{Step 4: Averaging transformation $\Phi_2$.}
	To reduce the order of the terms containing $t$ in the Hamiltonian \eqref{eq22}, introduce the canonical transformation $\Phi_{2}:(r_2,\theta_2,t)\rightarrow(r_1,\theta_1,t)$ as follows
	\begin{equation*}
		\left\{\begin{array}{l}
			r_1=r_2+\frac{\partial S_1}{\partial \theta_1}(r_2,\theta_1,t),  \\ [0.25cm]
			\theta_2=\theta_1+\frac{\partial S_1}{\partial r_2}(r_2,\theta_1,t), \end{array}
		\right.	
	\end{equation*}
	where
	\begin{equation*}
		S_1(r_2,\theta_1,t)=-\rho m\left(\frac{n}{m}\right)^\frac{3}{2}Q(\theta_1,t)r_2^\frac{1}{2},
	\end{equation*}
	and $Q$ is defined by \eqref{QQ}.
	
	It is worth noting that the assumption $\Phi_f(\theta)\equiv 0$ implies
	\begin{equation}\label{000}
		\int_{0}^{2\pi}C(ns)f(\theta_1+ns)ds\equiv 0,
	\end{equation}
	and consequently this transformation preserves the periodicity and boundedness of solutions. Indeed, $\Phi_f(\theta)\equiv 0$ gives
	\begin{equation*}
		\begin{aligned}
			&0\equiv\int_{0}^{2\pi}C\left(\frac{n}{m}\theta+nt\right)f(nt)dt=\int_{0}^{2\pi}C(ns)f\left(ns-\frac{n}{m}\theta\right)ds.
		\end{aligned}
	\end{equation*}
	Therefore, \eqref{000} holds.
	
	Under the transformation $\Phi_{2}$, the new Hamiltonian is
	\[
	H_2(r_2,\theta_2,t)
	=
	P_1
	+
	P_2r_2^{-\frac12}
	+
	P_3r_2^{-1}
	+
	O\left(r_2^{-\frac32}\right),
	\]
	
	where
	\begin{equation*}
		P_1(\theta_2,t)=\frac{n^2}{2m}\rho^2C^2(nt)f^2(\theta_2+nt)-\frac{\kappa^2}{2}C(nt)f(\theta_2+nt)\partial_{\theta_2}Q(\theta_2,t),
	\end{equation*}
	and $P_2$ and $P_3$ are given by \eqref{P_2} and \eqref{P_3}.
	
	\medskip
	\noindent\textbf{Step 5: Averaging transformation $\Phi_3$.}
	Next, we perform another canonical transformation $\Phi_3:(r_3,\theta_3,t)\mapsto(r_2,\theta_2,t)$ as follows
	\[\left
	\{\begin{array}{l}
		r_2=r_3+\frac{\partial S_2}{\partial \theta_2}(r_3,\theta_2,t),  \\ [0.25cm]
		\theta_3=\theta_2+ \frac{\partial S_2}{\partial r_3}(r_3,\theta_2,t),
	\end{array}
	\right.	\]
	where
	\[S_2(r_3,\theta_2,t)=\int_{0}^{t}\left([P_1](\theta_2)-P_1(\theta_2,s)\right)\,ds,\]
	\[[P_1](\theta_2)=\frac{1}{2\pi}\int_{0}^{2\pi}P_1(\theta_2,t)dt.\]
	Under the transformation $\Phi_3$, the Hamiltonian $H_2$ is transformed into
	\begin{equation*}
		\begin{aligned}
			H_3(r_3,\theta_3,t)&=H_2(r_2,\theta_2,t)+\frac{\partial}{\partial t}S_2(\theta_2,t)\\
			&=[P_1](\theta_3)+P_2(\theta_3,t)r_3^{-\frac12}
			+P_3(\theta_3,t)r_3^{-1}+O\left(r_3^{-\frac32}\right).
		\end{aligned}
	\end{equation*}

	\medskip
	\noindent\textbf{Step 6: Averaging transformation $\Phi_4$.}
	We first construct the canonical transformation $\Phi_4$ to eliminate
	the oscillatory part of $P_2(\theta_3,t)$ with respect
	to the time variable $t$.
	
	Recall that
	\[
	[P_2](\theta)
	=
	\frac{1}{2\pi}
	\int_0^{2\pi}
	P_2(\theta,t)\,dt
	\]
	and define
	\[
	W(\theta,t)
	=
	\int_0^t
	\left(
	[P_2](\theta)
	-
	P_2(\theta,s)
	\right)\,ds.
	\]
	Since the integrand has zero average with respect to $s$,
	the function $W(\theta,t)$ is $2\pi$-periodic in $t$. Choose the generating function of $\Phi_4$ as
	\[
	S_3(\theta_3,r_4,t)
	=W(\theta_3,t)r_4^{-\frac12}.
	\]
	Therefore, the canonical transformation $\Phi_4$ can be written as
	\[
	\Phi_4:
	\begin{cases}
		\displaystyle
		r_3
		=
		r_4+W_\theta(\theta_3,t)r_4^{-\frac12},\\[2mm]
		\displaystyle
		\theta_4
		=
		\theta_3-\frac12W(\theta_3,t)r_4^{-\frac32}.
	\end{cases}
	\]
	
	Under this time-dependent canonical transformation, the new Hamiltonian is
	\[
	H_4
	=
	H_3\circ\Phi_4
	+
	\frac{\partial}{\partial t}
	\left(
	W(\theta_3,t)r_4^{-\frac12}
	\right).
	\]
	Notice that
	\[
	\frac{\partial W}{\partial t}
	=
	[P_2](\theta_3)
	-
	P_2(\theta_3,t).
	\]
	Consequently,
	\[
	\begin{aligned}
		H_4(r_4,\theta_4,t)
		=[P_1](\theta_4)
		+
		[P_2](\theta_4)r_4^{-\frac12}+
		P_3(\theta_4,t)r_4^{-1}
		+
		O\left(r_4^{-\frac32}\right).
	\end{aligned}
	\]

	\medskip
	\noindent\textbf{Step 7: Averaging transformation $\Phi_5$.}
	We next construct the canonical transformation $\Phi_5$ to eliminate
	the oscillatory part of $P_3(\theta_4,t)$ with respect to
	the time variable $t$.
	
	Recall that
	\[
	[P_3](\theta)
	=
	\frac{1}{2\pi}
	\int_0^{2\pi}
	P_3(\theta,t)\,dt
	\]
	and define
	\[
	V(\theta,t)
	=
	\int_0^t
	\left(
	[P_3](\theta)
	-
	P_3(\theta,s)
	\right)\,ds.
	\]
	Similarly, $V(\theta,t)$ is $2\pi$-periodic in $t$. Choose the generating function of $\Phi_5$ as
	\[
	S_4(\theta_4,r_5,t)=V(\theta_4,t)r_5^{-1}.
	\]
	Then the canonical transformation generated by $\mathcal S_4$ can be written as
	\[
	\Phi_5:
	\begin{cases}
		\displaystyle
		r_4
		=
		r_5+V_\theta(\theta_4,t)r_5^{-1},\\[2mm]
		\displaystyle
		\theta_5
		=
		\theta_4-V(\theta_4,t)r_5^{-2}.
	\end{cases}
	\]
	
	The new Hamiltonian is given by
	\[
	H_5
	=
	H_4\circ\Phi_5
	+
	\frac{\partial}{\partial t}
	\left(
	V(\theta_4,t)r_5^{-1}
	\right).
	\]
	Since
	\[
	\frac{\partial V}{\partial t}
	=
	[P_3](\theta_4)
	-
	P_3(\theta_4,t),
	\]
	we obtain
	\[
	\begin{aligned}
		H_5(r_5,\theta_5,t)
		={}
		[P_1](\theta_5)
		+
		[P_2](\theta_5)r_5^{-\frac12}
		+
		[P_3](\theta_5)r_5^{-1}
		+
		O\left(r_5^{-\frac32}\right).
	\end{aligned}
	\]
	
	For notational simplicity, relabeling $(r_5,\theta_5,H_5)$ as $(r,\theta,H)$, we obtain the final normal form
	\begin{equation}\label{simple}
		\begin{aligned}
			H(r,\theta,t)
			=[P_1](\theta)+[P_2](\theta)r^{-\frac12}+[P_3](\theta)r^{-1}+O\left(r^{-\frac32}\right).
		\end{aligned}
	\end{equation}

	\section{Unbounded solutions under nonconstant $[P_1]$ }
	The main purpose of this section is to prove the existence of an unbounded solution for equation \eqref{equ} when $[P_1]$ is nonconstant. The discussion throughout this section is based on the assumptions that $f$ is a real-analytic $2\pi$-periodic function, $\Phi_{f}(\theta)\equiv0$, and $[P_1]$ is not constant. We first discuss the properties of $[P_1]$.
	
	First, since $f$ is real analytic and $2\pi$-periodic, $[P_1](\theta)$ is also a real-analytic function of $\theta$ that is $2\pi$-periodic, and thus $[P_1]'(\theta)$ is likewise periodic in $\theta$ with period $2\pi$. In particular, the mean of $[P_1]'(\theta)$ with respect to $\theta$ is zero, and hence $ [P_1]'(\theta)$ has at least one zero. Since $[P_1](\theta)$ is a nonconstant real-analytic periodic function, it attains a maximum at some point $\theta^*$. By the isolated-zero theorem for analytic functions, there exists an even integer $k\geq2$ such that, for $1\leq i\leq k-1$, $[P_1]^{(i)}(\theta^*)=0$, while $[P_1]^{(k)}(\theta^*)<0$. Consequently, $\partial_\theta[P_1](\theta^*+\tau)<0$ for all sufficiently small \(\tau>0\).
	
	The following lemma provides estimates for an analytic function near a zero.
	
	\begin{lemma}\label{lemma1}
		Suppose that $g(x)$ is analytic at $x=0$ and that $g^{(j)}(0)=0$, $j=0,1,\dots,k-1$, and $g^{(k)}(0)<0$. For any $c_1>1$ and $0<c_2<1$, there exists $\delta_0>0$ such that, whenever $0\leq x\leq \delta_0$, the following estimate holds:
		\[ g(x)\in \left[\frac{c_1}{k!}x^{k}g^{(k)}(0),\frac{c_2}{k!}x^{k}g^{(k)}(0)\right]. \]
	\end{lemma}
	\begin{proof}
		On the one hand, let \[ h_1(x)=g(x)-\frac{c_1}{k!}x^{k}g^{(k)}(0).\]
		It is then easy to obtain $h_1(0)=h_1^{\prime}(0)=\dots=h_1^{(k-1)}(0)=0$ and $h_1^{(k)}(0)=(1-c_1)g^{(k)}(0)>0$. Therefore, there exists $\delta_1>0$ such that, for $x\in [0,\delta_1]$, one has $h_1(x)\geq 0$. This implies \[ g(x)\geq\frac{c_1}{k!}x^{k}g^{(k)}(0).\]
		
		On the other hand, let \[ h_2(x)=g(x)-\frac{c_2}{k!}x^{k}g^{(k)}(0).\]
		Similarly, $h_2(0)=h_2^{\prime}(0)=\dots=h_2^{(k-1)}(0)=0$ and $h_2^{(k)}(0)=(1-c_2)g^{(k)}(0)<0$. Hence there exists $\delta_2>0$ such that, for $x\in [0,\delta_2]$, one has $h_2(x)\leq 0$. It follows that, for $x\in [0,\delta_2]$, \[ g(x)\leq\frac{c_2}{k!}x^{k}g^{(k)}(0).\]
		
		Now set $\delta_0:=\min\{\delta_1,\delta_2\}>0$. We immediately obtain the following conclusion:
		if $x$ satisfies $0\leq x\leq \delta_0$, then
		\[ g(x)\in \left[\frac{c_1}{k!}x^{k}g^{(k)}(0),\frac{c_2}{k!}x^{k}g^{(k)}(0)\right]. \]  \end{proof}
	
	Let $\tau=\theta-\theta^*$. Then the Hamiltonian system \eqref{simple} can be rewritten as
	\begin{equation}\label{eq18}
		\left\{\begin{array}{l}
			\frac{dr}{dt}=-\partial_{\tau}[P_1](\theta^*+\tau)-\partial_{\tau} R_1(r,\theta^*+\tau,t), \\ [0.25cm]
			\frac{d\tau}{dt}=\partial_r R_1(r,\theta^*+\tau,t), \end{array}
		\right.	
	\end{equation}
	where $R_1(r,\theta,t)=[P_2](\theta)r^{-\frac12}+[P_3](\theta)r^{-1}+O\left(r^{-\frac32}\right)$ in \eqref{simple}, and the remainder estimates are uniform and satisfy $\partial_\theta R_1(r,\theta,t)=O(r^{-\frac12})$ and $\partial_rR_1(r,\theta,t)=O(r^{-\frac32})$ for $r$ sufficiently large.
	
	The preceding discussion immediately yields Lemma \ref{lemma2} below.
	\begin{lemma}\label{lemma2}
		There exists $\delta>0$ such that, for every $\tau$ satisfying $0<\tau\le\delta$,
		\begin{equation*}
			\partial_\tau[P_1](\theta^*+\tau)\in\left[\frac{c_1}{(k-1)!}\tau^{k-1}[P_1]^{(k)}(\theta^*),\frac{c_2}{(k-1)!}\tau^{k-1}[P_1]^{(k)}(\theta^*)\right],
		\end{equation*}
		where $k\geq2$ is the positive integer satisfying $[P_1]'(\theta^*)=\cdots=[P_1]^{(k-1)}(\theta^*)=0$, $[P_1]^{(k)}(\theta^*)<0$, while $c_1>1$ and $0<c_2<1$ are taken from Lemma \ref{lemma1}.
	\end{lemma}
	\begin{proof}
		Let $g(\tau)=\partial_\theta[P_1](\theta^*+\tau)$. Then \[g^{(j)}(0)=[P_1]^{(j+1)}(\theta^*)=0,\quad j=0,1,\dots,k-2, \] and \[g^{(k-1)}(0)=[P_1]^{(k)}(\theta^*)<0.\] By Lemma \ref{lemma1}, we obtain
		\[ g(\tau)\in \left[\frac{c_1}{(k-1)!}\tau^{k-1}g^{(k-1)}(0),\frac{c_2}{(k-1)!}\tau^{k-1}g^{(k-1)}(0)\right],\] that is,
		\[ \partial_\theta[P_1](\theta^*+\tau)\in\left[\frac{c_1}{(k-1)!}\tau^{k-1}[P_1]^{(k)}(\theta^*),\frac{c_2}{(k-1)!}\tau^{k-1}[P_1]^{(k)}(\theta^*)\right]. \]
	\end{proof}
	
	In the sequel, we will prove the existence of unbounded solutions using the following iteration lemma.
	
	\begin{lemma}\label{iterate}
		There exist positive constants $r^*$, $C$ and $C_1$ such that if
		\[
		r(0)\geqslant r^*,\qquad r(0)^{-\frac{1}{2k+1}}\leqslant\tau(0)\leqslant \delta_0,
		\]
		then, for every $n=0,1,2,\dots$,
		\[
		r(2(n+1)\pi)^{\frac{3}{2}-\varepsilon}\geqslant r(2n\pi)^{\frac{3}{2}-\varepsilon}+2\pi C,
		\]
		and
		\[r(2(n+1)\pi)^{-\frac{1}{2k+1}}\leqslant\tau(2(n+1)\pi)\leqslant \delta_{n+1},\]
		where $\varepsilon=\frac{3}{2(2k+1)}$, $\delta_0=\frac{\delta}{2}$, $\delta_n=\delta_0+2\pi C_1 \dsum_{j=1}^{n}\left(\frac{r(2j\pi)}{2}\right)^{-\frac{3}{2}}$, $C$ and $C_1$ are two positive constants.
	\end{lemma}
	
	\begin{proof}
		We prove the result by induction. We first consider the case $n=0$, that is, $t\in[0,2\pi]$. It follows from \eqref{eq18} that, for $t\in[0,2\pi]$, $\frac{d\tau}{dt}=O(r^{-\frac{3}{2}}).$ Therefore, there exists a positive constant $C_1$ such that
		\begin{equation}\label{tau}
			\left|\frac{d\tau}{dt}\right|<C_1\left(\inf_{t\in[0,2\pi]}r(t)\right)^{-\frac{3}{2}}.
		\end{equation}
		Similarly, \eqref{eq18} shows that, for $t\in[0,2\pi]$, when $r(0)\geqslant r^*$ is sufficiently large, $\frac{dr}{dt}=O(1)$. Thus there exists a positive constant $C_2>0$ such that
		\begin{equation}\label{r}
			\left|\frac{dr}{dt}\right|<C_2.
		\end{equation}
		It follows from \eqref{r} that
		\[|r(t)-r(0)|\leqslant \left|\int_{0}^{t}\frac{dr}{ds}\,ds\right|\leqslant \int_{0}^{t} C_2\,ds\leqslant 2\pi C_2 ,\ t\in[0,2\pi].
		\]
		Hence, by taking $r(0)\geqslant r^*$ sufficiently large, we have, for every $t\in[0,2\pi]$,
		\begin{equation*}
			\frac{1}{2}r(0)<r(0)-2\pi C_2<r(t).
		\end{equation*}
		
		Moreover, by \eqref{tau}, for $t\in[0,2\pi]$, \[|\tau(t)-\tau(0)|\leqslant \int_{0}^{t}\left|\frac{d\tau}{ds}\right|\,ds\leqslant 2\pi C_1\left(\inf_{t\in[0,2\pi]} r(t)\right)^{-\frac{3}{2}}\leqslant 2\pi C_1\left(\frac{r(0)}{2}\right)^{-\frac{3}{2}}.\]
		Thus, for $t\in[0,2\pi]$,
		\begin{equation}\label{n_0}
			\begin{aligned}
				r(0)^{-\frac{1}{2k+1}}-2\pi C_1\left(\frac{r(0)}{2}\right)^{-\frac{3}{2}}\leqslant \tau(0)-2\pi C_1\left(\frac{r(0)}{2}\right)^{-\frac{3}{2}}\leqslant\tau(t),\\
				\tau(t)\leqslant \tau(0)+2\pi C_1\left(\frac{r(0)}{2}\right)^{-\frac{3}{2}}\leqslant \frac{\delta}{2}+2\pi C_1\left(\frac{r(0)}{2}\right)^{-\frac{3}{2}}.
			\end{aligned}
		\end{equation}
		Observe that, when $r^*$ is sufficiently large,
		\begin{equation*}
			\frac{1}{2}r(0)^{-\frac{1}{2k+1}}\leqslant r(0)^{-\frac{1}{2k+1}}-2\pi C_1\left(\frac{r(0)}{2}\right)^{-\frac{3}{2}} ,\ 2\pi C_1\left(\frac{r(0)}{2}\right)^{-\frac{3}{2}}\leqslant \frac{\delta}{2}.
		\end{equation*}
		Therefore, $r^*$ can be chosen sufficiently large so that, for $t\in[0,2\pi]$,
		\begin{equation*}
			0<\tau(t)\leqslant \delta.
		\end{equation*}
		
		We next prove that, for every $t\in[0,2\pi]$, $\tau(t)\geqslant r(t)^{-l}$, where $l=\frac{1}{2k+1}$. Suppose otherwise. Then there exists $\bar{t}\in(0,2\pi]$ such that $\tau(\bar{t})<r(\bar{t})^{-l}$. Consequently, there exists
		\begin{equation*}
			t_1=\sup\left\{s:\tau(s)\geqslant r(s)^{-l}, 0\le s<\bar{t}\right\}.
		\end{equation*}
		By the definition of $t_1$,
		\begin{equation}\label{relation}
			r(t_1)^{l}\tau(t_1)=1.
		\end{equation}
		Since $R_1(r,\theta,t)=O\left(r^{-\frac{1}{2}}\right)$, we have $r^{l}\partial_rR_1(r,\theta,t)=O\left(r^{l-\frac{3}{2}}\right)$. Hence, by Lemma \ref{lemma2} and relation \eqref{relation} satisfied by $t_1$,
		\begin{equation}\label{eq9}
			\begin{split}
				(r(t)^{l}\tau(t))^{'}|_{t=t_1}&=lr(t_1)^{l-1}r^{'}(t_1)\tau(t_1)+r(t_1)^{l}\tau^{'}(t_1)\\
				&=lr(t_1)^{l-1}\tau(t_1)(-\partial_{\tau}[P_1](\theta^*+\tau)-\partial_{\tau} R_1(r,\theta,t_1))+r(t_1)^{l}\partial_{r}R_1(r,\theta,t_1)\\
				&=lr(t_1)^{-1}(-\partial_{\tau}[P_1](\theta^*+\tau)-\partial_{\tau} R_1)+r(t_1)^{l}\partial_{r}R_1\\
				&\ge \frac{-lc_2}{(k-1)!}r(t_1)^{-1}\tau(t_1)^{k-1}[P_1]^{(k)}(\theta^*)+O\left(r^{l-\frac{3}{2}}\right)\\
				&= \frac{-lc_2}{(k-1)!}r(t_1)^{-[(k-1)l+1]}[P_1]^{(k)}(\theta^*)+O\left(r^{l-\frac{3}{2}}\right).
			\end{split}
		\end{equation}
		Observe that $-(k-1)l-1>l-\frac{3}{2}$. Therefore, when $r(0)\geqslant r^*$ is sufficiently large, the first term in \eqref{eq9} is dominant, and
		\[\frac{-lc_2}{(k-1)!}r(t_1)^{-[(k-1)l+1]}[P_1]^{(k)}(\theta^*)+O(r^{l-\frac{3}{2}})>0.\] Thus $(r(t)^{l}\tau(t))^{'}|_{t=t_1}>0$. Hence there exists $t_2>t_1$ such that $r(t_2)^l\tau(t_2)> 1$, which contradicts the definition of $t_1$. We have therefore proved that, for every $t\in[0,2\pi]$,
		\begin{equation}\label{2k+1}
			\tau(t)\geqslant r^{-\frac{1}{2k+1}}(t).
		\end{equation}
		
		We next prove that $\tau (2\pi)\leqslant \tau(0)+2\pi C_1\left(\frac{r(2\pi)}{2}\right)^{-\frac{3}{2}}$. By Lemma \ref{lemma2} and \eqref{2k+1}, for sufficiently large $r^*$,
		\begin{equation*}
			\begin{aligned}
				\frac{dr}{dt}&=-\partial_{\tau}[P_1](\theta^*+\tau)-\partial_{\tau} R_1(r,\theta^*+\tau,t)\\
				&\geqslant -\frac{c_2}{2(k-1)!}r^{-\frac{k-1}{2k+1}} [P_1]^{(k)}(\theta^*),\ t\in[0,2\pi].
			\end{aligned}
		\end{equation*}
		Write $-\frac{k-1}{2k+1}=-\frac{1}{2}+\varepsilon$, where $\varepsilon=\frac{3}{2(2k+1)}$.
		Consequently, there exists $C>0$ such that
		\begin{equation*}
			\frac{dr^{\frac{3}{2}-\varepsilon}}{dt}\geqslant -\frac{c_2}{2(k-1)!} [P_1]^{(k)}(\theta^*)\geqslant C.
		\end{equation*}
		Integrating the preceding inequality from $t=0$ to $t=2\pi$ gives
		\begin{equation*}
			r(2\pi)^{\frac{3}{2}-\varepsilon}-r(0)^{\frac{3}{2}-\varepsilon}\geqslant 2\pi C.
		\end{equation*}
		Clearly, $r(2\pi)>r(0)\geqslant r^*$. Integrating \eqref{r} for $t\in[0,2\pi]$ up to $2\pi$ gives
		\[
		|r(t)-r(2\pi)|\leqslant \left|\int_{t}^{2\pi}\frac{dr}{ds}\,ds\right|\leqslant \int_{t}^{2\pi} C_2\,ds\leqslant 2\pi C_2 ,\ t\in[0,2\pi].
		\]
		Therefore, since $r(2\pi)>r(0)$, we obtain
		\begin{equation*}
			\frac{r(2\pi)}{2}\leqslant r(t),\ t\in[0,2\pi].
		\end{equation*}
		Hence, by \eqref{tau},
		\begin{equation*}
			\tau(2\pi)\leqslant \tau(0)+2\pi C_1 \left(\frac{r(2\pi)}{2}\right)^{-\frac{3}{2}}\leqslant \delta_0+2\pi C_1 \left(\frac{r(2\pi)}{2}\right)^{-\frac{3}{2}}=\delta_1.
		\end{equation*}
		This completes the proof for the case $n=0$.
		
		Assume that the iteration lemma holds for \(0\le i\le n-1\), where \(n\ge1\), and consider $t\in[2n\pi,2(n+1)\pi]$. By the conclusion of the iteration lemma, $r(2n\pi)>\cdots>r(2\pi)>r(0)\geqslant r^*$ is a strictly increasing sequence, and
		\begin{equation*}
			\begin{aligned}
				\delta_n&=\delta_0+2\pi C_1 \sum_{j=1}^{n}\left(\frac{r(2j\pi)}{2}\right)^{-\frac{3}{2}}\\
				&=\delta_0+2^{\frac{3}{2}}2\pi C_1 \sum_{j=1}^{n}\left(r(2j\pi)^{\frac{3}{2}-\varepsilon}\right)^{-\frac{1}{1-\frac{2}{3}\varepsilon}}\\
				&\leqslant \delta_0+2^{\frac{3}{2}}2\pi C_1 \sum_{j=1}^{n}\left(r(0)^{\frac{3}{2}-\varepsilon}+2j\pi C\right)^{-\frac{1}{1-\frac{2}{3}\varepsilon}}:=\delta_0+S_n.
			\end{aligned}
		\end{equation*}
		Therefore, $r^*$ can be chosen sufficiently large so that $S_\infty\leqslant\frac{\delta}{4}$.
		
		Since $r(2n\pi)>\cdots>r(2\pi)>r(0)\geqslant r^*$, estimate \eqref{n_0} also holds for $i=n$; namely, for $t\in[2n\pi,2(n+1)\pi]$,
		\begin{equation*}
			0<r(2n\pi)^{-\frac{1}{2k+1}}-2\pi C_1\left(\frac{r(2n\pi)}{2}\right)^{-\frac{3}{2}}\leqslant \tau(2n\pi)-2\pi C_1\left(\frac{r(2n\pi)}{2}\right)^{-\frac{3}{2}}\leqslant\tau(t),
		\end{equation*}
		and
		\begin{equation*}
			\tau(t)\leqslant\tau(2n\pi)+2\pi C_1\left(\frac{r(2n\pi)}{2}\right)^{-\frac{3}{2}}\leqslant \delta_n+2\pi C_1\left(\frac{r(2n\pi)}{2}\right)^{-\frac{3}{2}}\leqslant \delta_0+2S_n\leqslant \delta .
		\end{equation*}
		
		Therefore, by Lemma \ref{lemma2} and an argument similar to that for $n=0$, one can prove that, for every $t\in[2n\pi,2(n+1)\pi]$, $\tau(t)\geqslant r^{-\frac{1}{2k+1}}(t)$, and that
		\begin{equation*}
			r(2(n+1)\pi)^{\frac{3}{2}-\varepsilon}-r(2n\pi)^{\frac{3}{2}-\varepsilon}\geqslant 2\pi C,
		\end{equation*}
		as well as $\tau(2(n+1)\pi)\leqslant \delta_{n+1}$.
	\end{proof}
	
	\medskip
	\noindent\textbf{Proof of part one of Theorem \ref{th}.}
	Since $f$ is a real-analytic $2\pi$-periodic function satisfying $\Phi_{f}(\theta)\equiv0$, and $[P_1]$ is not constant, Lemma \ref{lemma2} applies. Take $r(0)\geqslant r^*$ sufficiently large. By Lemma \ref{iterate}, if $(r(0),\tau(0))\in \{(r,\tau):r\geq r^*,r^{-\frac{1}{2k+1}}\leq\tau\leq\delta_0\}$, then, for every $n\geqslant 0$,
	\[
	r(2(n+1)\pi)^{\frac{3}{2}-\varepsilon}\geqslant r(2n\pi)^{\frac{3}{2}-\varepsilon}+2\pi C.
	\]
	Iterating this inequality gives
	\[
	r(2n\pi)^{\frac{3}{2}-\varepsilon}\geqslant r(0)^{\frac{3}{2}-\varepsilon}+2n\pi C\longrightarrow\infty.
	\]
	Therefore, the solution $(r(t), \tau(t))$ is unbounded. This completes the proof of part one of Theorem \ref{th}.
	\hfill\qedsymbol

	\section{Unbounded solutions under degenerate $[P_1]$}
	In this section, we consider the case where both $\Phi_{f}(\theta)\equiv0$ and $[P_1](\theta)\equiv constant$. Under this assumption, the Hamiltonian $H$ in \eqref{simple} degenerates to a constant plus lower-order terms. Since an additive constant does not affect the Hamiltonian vector field, we write
	\begin{equation}\label{eq:normal-form}
		\begin{aligned}
			H(r,\theta,t)
			=[P_2](\theta)r^{-\frac12}+[P_3](\theta)r^{-1}+R(r,\theta,t).
		\end{aligned}
	\end{equation}
	Here the term $R$ satisfies the following property: for some sufficiently large integer \(s\), and for all nonnegative integers $i,j,\ell$ satisfying $i+j+\ell\leq s$, there is a constant $C_{ij\ell}>0$ such that
	\begin{equation}
		\label{eq:symbol-estimate}
		\left|
		\partial_r^i\partial_\theta^j\partial_t^\ell
		R(r,\theta,t)
		\right|
		\leq C_{ij\ell}r^{-3/2-i}
	\end{equation}
	uniformly for $(\theta,t)\in\mathbb{T}^2$ and all sufficiently large $r$.
	
	We now state the main result of this section. Theorem \ref{thm:main} directly yields part (2) of Theorem~\ref{th}.
	
	\begin{theorem}\label{thm:main}
		Suppose that there exists $\theta_*\in\mathbb{T}$ such that
		\begin{equation*}
			[P_2](\theta_*)=0,
			\qquad
			\alpha:=[P_2]'(\theta_*)\neq0,
		\end{equation*}
		then the Hamiltonian system associated with \eqref{eq:normal-form}
		has an unbounded solution.  More precisely, if $\alpha<0$, there is a forward orbit for which
			$r(t)\to+\infty$ as $t\to+\infty$; and if $\alpha>0$, there is a backward orbit for which
			$r(t)\to+\infty$ as $t\to-\infty$. In the forward case, let
		$(\theta_n,r_n):=(\theta(2\pi n),r(2\pi n))$, then
		\begin{align}
			\label{eq:theta-asymptotic}
			\theta_n-\theta_*
			&=-\frac{[P_3](\theta_*)}{[P_2]'(\theta_*)}r_n^{-1/2}
			+O(r_n^{-1}),\\
			\label{eq:r-asymptotic}
			r_n
			&\sim\bigl(-3\pi [P_2]'(\theta_*)n\bigr)^{2/3}.
		\end{align}
		Consequently, whenever \eqref{eq:normal-form} is obtained from the asymmetric oscillator by canonical transformations asymptotic to the identity, the original equation possesses an unbounded
		solution.
	\end{theorem}

	The main point is the construction of orbits that lie on the one-dimensional stable manifold
	of the saddle point. We explain why, after the degeneracy of $\bigl[P_1\bigr]$,
	the argument used in the nondegenerate case cannot be applied directly.
	
	For simplicity, denote
	\[
	A(\theta)=[P_2](\theta),
	\qquad
	B(\theta)=[P_3](\theta).
	\]
	Introduce
	\[
	u=r^{-1/2}.
	\]
	Thus, $u=0$ corresponds to $r=+\infty$. The Hamiltonian equations yield
	
	\begin{equation}\label{u_theta_t_hamilton}
		\begin{cases}
			\displaystyle
			\frac{du}{dt}
			=\frac12 A'(\theta)u^4
			+\frac12 B'(\theta)u^5
			+O(u^6),
			\\[2mm]
			\displaystyle
			\frac{d\theta}{dt}=
			-\frac12 A(\theta)u^3
			-B(\theta)u^4
			+O(u^5).
		\end{cases}
	\end{equation}
	For $u>0$, we introduce the new time variable $s$ by
	\[
	\frac{ds}{dt}=u^3.
	\]
	Since this time transformation preserves the orientation of the
	orbits, the transformed system becomes
	\begin{equation}
		\begin{cases}
			\displaystyle
			\frac{du}{ds}
			=
			\frac12 A'(\theta)u
			+\frac12 B'(\theta)u^2
			+O(u^3),
			\\[2mm]
			\displaystyle
			\frac{d\theta}{ds}
			=
			-\frac12 A(\theta)
			-B(\theta)u
			+O(u^2).
		\end{cases}
		\label{eq:desingularized-system}
	\end{equation}
	
	Recall that \(\theta_*\) is a nondegenerate zero of $A$, namely,
	\[
	A(\theta_*)=0,
	\qquad
	\alpha:=A'(\theta_*)\neq0,
	\]
	and put
	\[
	q=B(\theta_*),
	\qquad
	\tau=\theta-\theta_*.
	\]
	Then \eqref{eq:desingularized-system} near $(0,0)$ can be written as
	\begin{equation*}
		\frac{d}{ds}
		\begin{pmatrix}
			u\\ \tau
		\end{pmatrix}
		=
		\begin{pmatrix}
			\dfrac{\alpha}{2} & 0\\[1mm]
			-q & -\dfrac{\alpha}{2}
		\end{pmatrix}
		\begin{pmatrix}
			u\\ \tau
		\end{pmatrix}
		+
		O\bigl((|u|+|\tau|)^2\bigr).
	\end{equation*}
	The eigenvalues of the linear part are
	\[
	\lambda_1=\frac{\alpha}{2},
	\qquad
	\lambda_2=-\frac{\alpha}{2},
	\]
	and they have opposite signs. Consequently, the equilibrium
	$(u,\tau)=(0,0)$ is a hyperbolic saddle point.

	To see the obstruction more explicitly, introduce the coordinate
	\[
	z=\tau+\frac{q}{\alpha}u.
	\]
	This change of coordinates diagonalizes the linear part, and the
	system becomes
	\begin{equation*}
		\begin{cases}
			\displaystyle
			\frac{du}{ds}
			=
			\frac{\alpha}{2}u+F(u,z,t(s)),
			\\[2mm]
			\displaystyle
			\frac{dz}{ds}
			=
			-\frac{\alpha}{2}z+G(u,z,t(s)),
		\end{cases}
	\end{equation*}
	where
	\[
	F(u,z,t(s)),\,G(u,z,t(s))
	=
	O\bigl((|u|+|z|)^2\bigr).
	\]
	The preceding estimates are uniform in $t\in\mathbb T$. Suppose that $\alpha<0$. Then the $u$-direction is stable, whereas the $z$-direction
	is unstable. This saddle-point structure explains why the invariant-region
	argument used before the degeneracy of $[P_1]$ cannot be
	copied. Hence only initial data lying on the one-dimensional local stable manifold can satisfy $z\rightarrow 0$ and $u\rightarrow0$; generic nearby orbits possess a nonzero unstable component and leave the neighborhood. Therefore, after the degeneracy of $[P_1]$, the existence
	of an unbounded solution is no longer obtained from a
	two-dimensional invariant region. Instead, it follows from the
	stable manifold of the saddle point created by the desingularizing
	time transformation. If $\alpha>0$, the same argument applied
	in backward time constructs the corresponding unstable manifold and
	an orbit unbounded in negative time.
	
	However, the classical continuous-time stable manifold theory cannot be applied directly, because the time rescaling makes the extended vector field
	\begin{equation*}
		\begin{cases}
			\displaystyle
			\frac{du}{ds}
			=\frac12 A'(\theta)u
			+\frac12 B'(\theta)u^2
			+O(u^3),
			\\[2mm]
			\displaystyle
			\frac{d\theta}{ds}=
			-\frac12 A(\theta)
			-B(\theta)u
			+O(u^2),
			\\[2mm]
			\displaystyle
			\frac{dt}{ds}=u^{-3},
		\end{cases}
	\end{equation*}
	whose $t$-component is singular at $u=0$. We therefore consider the Poincar\'{e} map, which eliminates the periodic time variable while preserving the full dynamics. Since the resulting fixed point is parabolic, a discrete Lyapunov--Perron argument in a polynomially weighted sequence space is used to construct the desired unbounded orbit rigorously.
	
	We first give the required discrete Lyapunov--Perron lemma. Let $\kappa>0$, and consider a $C^2$ map
	\begin{equation}
		\label{eq:abstract-map}
		\mathcal P(u,w)
		=\left(
		u-\kappa u^4+F(u,w),\,
		(1+2\kappa u^3)w+G(u,w)
		\right).
	\end{equation}
	Assume that there are constants $C>0$, $\varepsilon>0$, and
	$\delta_0>0$ such that, whenever
	\begin{equation}
		\label{eq:local-domain}
		0<u\leq\varepsilon,
		\qquad |w|\leq\delta_0,
	\end{equation}
	one has
	\begin{align}
		\label{eq:F-estimates}
		|F(u,w)|&\leq Cu^5,
		&|\partial_uF(u,w)|&\leq Cu^4,
		&|\partial_wF(u,w)|&\leq Cu^5,\\
		\label{eq:G-estimates}
		|G(u,w)|&\leq Cu^4,
		&|\partial_uG(u,w)|&\leq Cu^3,
		&|\partial_wG(u,w)|&\leq Cu^4.
	\end{align}
	
	\begin{lemma}
		\label{lem:LP}
		Under \eqref{eq:abstract-map}--\eqref{eq:G-estimates}, for every
		sufficiently small $u_0>0$ there exists $w_0=O(u_0)$ such that the
		forward orbit
		\[
		(u_n,w_n)=\mathcal P^n(u_0,w_0)
		\]
		is defined for every $n\geq0$ and satisfies
		\begin{equation*}
			u_n\to0,
			\qquad
			w_n=O(u_n),
			\qquad
			u_n\sim(3\kappa n)^{-1/3}.
		\end{equation*}
	\end{lemma}
	
	\begin{proof}
		Fix a small $u_0>0$ and define the reference sequence
		\begin{equation*}
			q_n=(u_0^{-3}+3\kappa n)^{-1/3},
			\qquad n\geq0.
		\end{equation*}
		Thus $q_0=u_0$ and $q_n\sim(3\kappa n)^{-1/3}$ as $n\rightarrow\infty$.
		
		\medskip
		\noindent\textbf{Step 1: The weighted sequence space.}
		Let
		\[
		X_{u_0}
		=\left\{w=(w_n)_{n\geq0}:
		\|{w}\|_*:=\sup_{n\geq0}\frac{|w_n|}{q_n}<\infty
		\right\}.
		\]
		For a constant $M>0$, to be fixed below, write
		\[
		\mathcal B_M=\{w\in X_{u_0}:\|w\|_*\leq M\}.
		\]
		Given $w\in\mathcal B_M$, define $u[w]=(u_n[w])_{n\geq0}$ by
		\begin{equation}
			\label{eq:u-recursion}
			u_{n+1}=u_n-\kappa u_n^4+F(u_n,w_n),
			\qquad u_0\text{ fixed}.
		\end{equation}
		We shall choose $u_0$ small enough that $u_0\leq\varepsilon$ and $Mu_0<\delta_0$.  The constants
		below are uniform for $w\in\mathcal B_M$.
		
		We first estimate the radial sequence $(u_n)_{n\geq 0}$. If $u_0$ is sufficiently small, \eqref{eq:u-recursion} gives
		\begin{equation}
			\label{eq:u-q-comparison}
			c_0q_n\leq u_n\leq C_0q_n,
			\qquad n\geq0,
		\end{equation}
		where $c_0,C_0>0$ do not depend on $n$ or on $w\in\mathcal B_M$.
		Indeed, \eqref{eq:F-estimates} gives
		\[
		u_n-\frac{3\kappa}{2}u_n^4
		\leq u_{n+1}
		\leq u_n-\frac{\kappa}{2}u_n^4.
		\]
		Consequently $u_n$ remains positive and decreases. The mean value theorem then yields two positive constants
		$c,C$ such that
		\[
		c\leq u_{n+1}^{-3}-u_n^{-3}\leq C.
		\]
		After summation, this is exactly \eqref{eq:u-q-comparison}.  Moreover,
		$|w_n|\leq Mq_n\leq Mu_0<\delta_0$, so every evaluation remains in
		the domain \eqref{eq:local-domain}.
		
		We need a sharper comparison.  From \eqref{eq:u-recursion},
		\begin{equation}
			\label{eq:inverse-cube-increment}
			u_{n+1}^{-3}-u_n^{-3}=3\kappa+O(u_n).
		\end{equation}
		Using \eqref{eq:u-q-comparison} and summing gives
		\begin{align}
			u_n^{-3}
			&=u_0^{-3}+3\kappa n
			+O\left(\sum_{k=0}^{n-1}q_k\right)\notag\\
			\label{eq:sharp-radial}
			&=u_0^{-3}+3\kappa n
			+O\left((u_0^{-3}+3\kappa n)^{2/3}\right).
		\end{align}
		Indeed, the last estimate follows by integral comparison.  Therefore
		\begin{equation}
			\label{eq:u3-q3}
			u_n^3=q_n^3\bigl(1+O(q_n)\bigr)
			=q_n^3+O(q_n^4).
		\end{equation}
		All constants in \eqref{eq:sharp-radial}--\eqref{eq:u3-q3} are
		uniform on $\mathcal B_M$.
		
		\medskip
		\noindent\textbf{Step 2: The estimate of the kernel.}
		Set
		\begin{equation*}
			A_n[w]=1+2\kappa u_n[w]^3,
			\qquad
			K_{n,j}[w]=\prod_{\ell=n}^{j}A_\ell[w]^{-1},
			\quad j\geq n.
		\end{equation*}
		For $u_0$ small, $A_n>1$, we claim that
		\begin{equation}
			\label{eq:kernel-estimate}
			0<K_{n,j}[w]
			\leq C_K\left(\frac{q_{j+1}}{q_n}\right)^2,
			\qquad j\geq n,
		\end{equation}
		uniformly for $w\in\mathcal B_M$.
		
		To prove this, take logarithms.  By \eqref{eq:u3-q3},
		\begin{align*}
			\log K_{n,j}
			&=-\sum_{\ell=n}^{j}\log(1+2\kappa u_\ell^3)\\
			&=-2\kappa\sum_{\ell=n}^{j}q_\ell^3
			+O\left(\sum_{\ell=n}^{j}q_\ell^4\right)
			+O\left(\sum_{\ell=n}^{j}q_\ell^6\right).
		\end{align*}
		The two error sums are bounded uniformly in $j$.  Moreover, by integral comparison, we have
		\[
		\sum_{\ell=n}^{j}q_\ell^3
		=\sum_{\ell=n}^{j}\frac{1}{u_0^{-3}+3\kappa\ell}
		\geq\frac{1}{3\kappa}
		\log\frac{u_0^{-3}+3\kappa(j+1)}{u_0^{-3}+3\kappa n}.
		\]
		It follows that
		\[
		K_{n,j}
		\leq C_K
		\left(
		\frac{u_0^{-3}+3\kappa n}{u_0^{-3}+3\kappa(j+1)}
		\right)^{2/3},
		\]
		which is \eqref{eq:kernel-estimate}.
		
		\medskip
		\noindent\textbf{Step 3: Definition and boundedness of the Lyapunov--Perron operator.}
		The inverse of $$w_{n+1}=A_{n}w_n+G(u_n,w_n)$$ is $$w_n=A_n^{-1}w_{n+1}-A_n^{-1}G(u_n,w_n).$$ Therefore,
		\begin{equation*}
			w_{n}=K_{n,N}w_{N+1}-\sum_{j=n}^{N}K_{n,j}G(u_j,w_j).
		\end{equation*}
		Since every $w\in\mathcal B_M$ satisfies $w_N\to0$ and $K_{n,N}$ is bounded, letting $N\to\infty$ motivates the definition
		\begin{equation}
			\label{eq:LP-operator}
			(\mathcal Tw)_n
			=-\sum_{j=n}^{\infty}
			K_{n,j}[w]G(u_j[w],w_j)
		\end{equation}
		for $w\in\mathcal B_M$.
		
		By \eqref{eq:G-estimates}, \eqref{eq:u-q-comparison}, and
		\eqref{eq:kernel-estimate},
		\begin{align*}
			|(\mathcal Tw)_n|
			&\leq Cq_n^{-2}\sum_{j=n}^{\infty}q_{j+1}^2q_j^4 \leq Cq_n^{-2}\sum_{j=n}^{\infty}q_j^6\leq C_1q_n.
		\end{align*}
		The final inequality follows from
		\[
		\sum_{j=n}^{\infty}q_j^6
		\leq C(u_0^{-3}+3\kappa n)^{-1}=Cq_n^3.
		\]
		Thus the series in \eqref{eq:LP-operator} converges absolutely.  If
		$M\geq2C_1$, then
		\begin{equation*}
			\mathcal T(\mathcal B_M)\subseteq\mathcal B_M.
		\end{equation*}
		
		\medskip
		\noindent\textbf{Step 4: Lipschitz dependence of the radial
			sequence.}
		Let $w,\widetilde w\in\mathcal B_M$, and use tildes for all objects
		constructed from $\widetilde w$.  Put
		$\delta=\|{w-\widetilde w}\|_*$.  Subtracting the two radial
		recursions and using the mean-value theorem,
		\eqref{eq:F-estimates}, and \eqref{eq:u3-q3}, one obtains
		\begin{equation}
			\label{eq:radial-difference-recursion}
			u_{n+1}-\widetilde u_{n+1}
			=\bigl(1-4\kappa q_n^3+O(q_n^4)\bigr)
			(u_n-\widetilde u_n)+e_n,
			\qquad |e_n|\leq Cq_n^6\delta,
		\end{equation}
		with $u_0-\widetilde u_0=0$.  One has
		\begin{equation*}
			\prod_{\ell=k+1}^{n-1}
			\bigl|1-4\kappa q_\ell^3+O(q_\ell^4)\bigr|
			\leq C\left(\frac{q_n}{q_{k+1}}\right)^4.
		\end{equation*}
		Indeed, taking logarithms reduces the main term to
		$-4\kappa\sum q_\ell^3$, which equals
		$-\frac43\log\!\bigl((u_0^{-3}+3\kappa n)/
		(u_0^{-3}+3\kappa(k+1))\bigr)$ up to a uniformly bounded error.
		\eqref{eq:radial-difference-recursion} now gives
		\begin{align}
			|u_n-\widetilde u_n|
			&\leq Cq_n^4\sum_{k=0}^{n-1}q_k^2\,\delta
			\leq Cq_n^3\delta.
			\label{eq:radial-Lipschitz}
		\end{align}
		We used
		$\sum_{k=0}^{n-1}q_k^2\leq Cq_n^{-1}$.
		
		\medskip
		\noindent\textbf{Step 5: Lipschitz dependence of the kernel.}
		From \eqref{eq:radial-Lipschitz},
		\[
		|A_\ell-\widetilde A_\ell|
		\leq Cq_\ell^2|u_\ell-\widetilde u_\ell|
		\leq Cq_\ell^5\delta.
		\]
		Define
		\[
		L_{n,j}[w]
		:=
		\sum_{\ell=n}^{j}
		\log A_\ell[w].
		\]
		When \(u_0\rightarrow 0\), \(A_\ell[w]\) and
		\(A_\ell[\widetilde w]\) tend to  \(1\). Hence, by the mean-value theorem,
		\[
		\left|
		\log A_\ell[w]
		-
		\log A_\ell[\widetilde w]
		\right|
		\leq
		C
		\left|
		A_\ell[w]-A_\ell[\widetilde w]
		\right|.
		\]
		Again using the mean-value theorem together with
		$\sum_{\ell=n}^{\infty}q_\ell^5\leq Cq_n^2$ and \eqref{eq:kernel-estimate} yields
		\begin{equation}
			\label{eq:kernel-Lipschitz}
			\begin{aligned}
				\left|
				K_{n,j}[w]-K_{n,j}[\widetilde w]
				\right|
				&=
				\left|
				e^{-L_{n,j}[w]}
				-
				e^{-L_{n,j}[\widetilde w]}
				\right|\\
				&\leq
				\max\left\{
				K_{n,j}[w],
				K_{n,j}[\widetilde w]
				\right\}
				\left|
				L_{n,j}[w]-L_{n,j}[\widetilde w]
				\right|\leq C
				\left(\frac{q_{j+1}}{q_n}\right)^2q_n^2\delta.
			\end{aligned}
		\end{equation}

		\medskip
		\noindent\textbf{Step 6: Contraction.}
		By \eqref{eq:G-estimates} and \eqref{eq:radial-Lipschitz},
		\begin{align}
			|G(u_j,w_j)-G(\widetilde u_j,\widetilde w_j)|
			&\leq Cq_j^3|u_j-\widetilde u_j|
			+Cq_j^4|w_j-\widetilde w_j|\notag\\
			&\leq Cq_j^5\delta.
			\label{eq:G-Lipschitz-sequence}
		\end{align}
		Splitting
		$\mathcal Tw-\mathcal T\widetilde w$ into the difference of the
		kernels and the difference of the $G$ terms, and applying
		\eqref{eq:kernel-estimate}, \eqref{eq:kernel-Lipschitz}, and
		\eqref{eq:G-Lipschitz-sequence}, we obtain
		\begin{align*}
			|(\mathcal Tw)_n-(\mathcal T\widetilde w)_n|
			&\leq Cq_n^{-2}\sum_{j=n}^{\infty}q_j^7\,\delta
			+Cq_n^2q_n^{-2}\sum_{j=n}^{\infty}q_j^6\,\delta\\
			&\leq C(q_n^2+q_n^3)\delta.
		\end{align*}
		After division by $q_n$ and taking the supremum in $n$,
		\begin{equation*}
			\|{\mathcal Tw-\mathcal T\widetilde w}\|_*
			\leq Cu_0\|{w-\widetilde w}\|_*.
		\end{equation*}
		Shrinking $u_0$ so that $Cu_0<1$, we conclude that $\mathcal T$ is
		a strict contraction of $\mathcal B_M$. Since $\mathcal B_M$ is a closed subset of the Banach space $X_{u_0}$, it is complete.
		
		By the Banach fixed-point theorem there is a unique
		$w^*\in\mathcal B_M$ satisfying $\mathcal Tw^*=w^*$.  Let
		$u^*=u[w^*]$.  The fixed-point identity is
		\begin{equation}
			\label{eq:LP-fixed-point}
			w_n^*=-\sum_{j=n}^{\infty}
			K_{n,j}[w^*]G(u_j^*,w_j^*).
		\end{equation}
		Multiplying \eqref{eq:LP-fixed-point} by
		$A_n=1+2\kappa(u_n^*)^3$ and separating the term $j=n$ gives
		\[
		w_{n+1}^*=A_nw_n^*+G(u_n^*,w_n^*).
		\]
		Together with \eqref{eq:u-recursion}, this proves
		\[
		(u_{n+1}^*,w_{n+1}^*)=\mathcal P(u_n^*,w_n^*)
		\]
		for every $n\geq0$.  Moreover,
		$|w_n^*|\leq Mq_n=O(u_n^*)$.
		
		Finally, \eqref{eq:inverse-cube-increment} and
		\eqref{eq:u-q-comparison} imply
		\[
		(u_n^*)^{-3}
		=u_0^{-3}+3\kappa n+O(n^{2/3})
		=3\kappa n+o(n).
		\]
		Hence $u_n^*\sim(3\kappa n)^{-1/3}$, completing the proof.
	\end{proof}
	
	\medskip
	\noindent\textbf{Proof of Theorem~\ref{thm:main}.}
	
	\medskip
	\noindent\textbf{Step 1: Variable transformation.}
	Set
	\begin{equation*}
		u=r^{-1/2},
		\qquad
		\tau=\theta-\theta_*.
	\end{equation*}
	After the transformation, \eqref{u_theta_t_hamilton} becomes
	\begin{equation*}
		\begin{aligned}
			\dot u
			&=\frac12A'(\theta_*+\tau)u^4
			+\frac12B'(\theta_*+\tau)u^5+O(u^6),\\
			\dot\tau
			&=-\frac12A(\theta_*+\tau)u^3
			-B(\theta_*+\tau)u^4+O(u^5).
		\end{aligned}
	\end{equation*}
	The remainders and their derivatives required below are uniform in
	$t\in\mathbb{T}$.
	
	Write $B_*=B(\theta_*)$.  Since
	\[
	A(\theta_*+\tau)=\alpha\tau+O(\tau^2),
	\qquad
	A'(\theta_*+\tau)=\alpha+O(\tau),
	\]
	and $B(\theta_*+\tau)=B_*+O(\tau)$, in every region
	$|\tau|\leq Cu$ the system becomes
	\begin{equation}
		\label{eq:u-tau-local}
		\begin{aligned}
			\dot u&=\frac{\alpha}{2}u^4+O(u^5),\\
			\dot\tau&=-\frac{\alpha}{2}\tau u^3-B_*u^4+O(u^5).
		\end{aligned}
	\end{equation}
	
	\medskip
	\noindent\textbf{Step 2: The time-$T$ map.}
	Let $T=2\pi$ and denote by $\mathcal P$ the time-$T$ map on the
	section $t=0\pmod T$. Integrating \eqref{eq:u-tau-local} over the bounded interval $[0,T]$ gives
	\begin{equation}
		\label{eq:P-u-tau}
		\begin{aligned}
			u_1
			&=u+\frac{T\alpha}{2}u^4+O(u^5),\\
			\tau_1
			&=\tau-\frac{T\alpha}{2}\tau u^3-TB_*u^4+O(u^5).
		\end{aligned}
	\end{equation}
	These expansions hold in the \(C^1\) sense, uniformly in every cone $|\tau|\leq Cu$.
	
	Introduce the variables $v$ and $w$ as
	\begin{equation*}
		v=\frac{\tau}{u},
		\qquad
		v_*=-\frac{B_*}{\alpha},
		\qquad
		w=v-v_*.
	\end{equation*}
	Dividing the second equation of \eqref{eq:P-u-tau} by the first and
	using $\tau=uv$ gives
	\begin{equation}
		\label{eq:P-u-v}
		\begin{aligned}
			u_1
			&=u+\frac{T\alpha}{2}u^4+O(u^5),\\
			v_1
			&=v-T(\alpha v+B_*)u^3+O(u^4).
		\end{aligned}
	\end{equation}
	
	We first consider the case $\alpha<0$. Set
	\begin{equation*}
		\kappa=-\frac{T\alpha}{2}>0.
	\end{equation*}
	In the $(u,w)$ variables, \eqref{eq:P-u-v} has the form
	\begin{equation*}
		\begin{aligned}
			u_1&=u-\kappa u^4+F(u,w),\\
			w_1&=(1+2\kappa u^3)w+G(u,w).
		\end{aligned}
	\end{equation*}
	\eqref{eq:symbol-estimate} yields precisely
	\eqref{eq:F-estimates}--\eqref{eq:G-estimates} in every fixed bounded
	strip $|w|\leq\delta_0$.  Thus Lemma~\ref{lem:LP} applies.
	
	\medskip
	\noindent\textbf{Step 3: The escaping orbit.}
	For every sufficiently small $u_0>0$, Lemma~\ref{lem:LP} selects an
	initial value $w_0=O(u_0)$ whose orbit satisfies
	\begin{equation}
		\label{eq:selected-orbit}
		u_n\to0,
		\qquad
		w_n=O(u_n),
		\qquad
		u_n\sim(3\kappa n)^{-1/3}.
	\end{equation}
	Since $r_n=u_n^{-2}$, it follows that $r_n\to+\infty$.  Furthermore,
	\[
	v_n=v_*+O(u_n),
	\]
	and hence
	\[
	\theta_n-\theta_*
	=\tau_n=u_nv_n
	=-\frac{B_*}{\alpha}u_n+O(u_n^2).
	\]
	This proves \eqref{eq:theta-asymptotic}.  Since
	$\kappa=-T\alpha/2=-\pi\alpha$, \eqref{eq:selected-orbit} also gives
	\[
	r_n=u_n^{-2}
	\sim(3\kappa n)^{2/3}
	=(-3\pi\alpha n)^{2/3},
	\]
	which is \eqref{eq:r-asymptotic}.
	
	The estimates are uniform during each interval
	$[nT,(n+1)T]$.  In fact, \eqref{eq:u-tau-local} gives
	\[
	\sup_{0\leq s\leq T}|u(nT+s)-u_n|=O(u_n^4).
	\]
	Therefore $u(t)\to0$ and $r(t)\to+\infty$ as $t\to+\infty$.
	
	If $\alpha>0$, apply the preceding construction to the inverse
	time-$T$ map.  Its radial coefficient has the opposite sign, so the
	same Lyapunov--Perron argument constructs an orbit converging to
	$u=0$ under backward iteration.  Equivalently,
	$r(t)\to+\infty$ as $t\to-\infty$.
	
	\hfill\qedsymbol

	\section{Unbounded solutions under degenerate $[P_1]$ and $[P_2]$}
	
	In this section, suppose that
	\[
	[P_1](\theta)\equiv constant,\quad[P_2](\theta)\equiv 0.
	\]
	Since an additive constant does not affect the Hamiltonian vector field, we write the Hamiltonian as
	\begin{equation}\label{eq:hamiltonian-P3}
		H(r,\theta,t)
		=[P_3](\theta)r^{-1}
		+[P_4](\theta)r^{-3/2}
		+O(r^{-2}),
	\end{equation}
	where $[P_4]$ is the coefficient of the next-order term obtained by performing a further canonical transformation and the remainder and its first derivatives satisfy the corresponding uniform estimates.

	We show that if $[P_3]$ has a nondegenerate
	zero, then the preceding Lyapunov--Perron argument can still be applied. The following theorem directly yields part (3) of Theorem~\ref{th}.

	\begin{theorem}
		\label{thm:P3-nondegenerate-zero}
		Suppose that there
		exists $\theta_*\in \mathbb{T}$ such that
		\[
		\bigl[P_3\bigr](\theta_*)=0,
		\qquad
		\bigl[P_3\bigr]'(\theta_*)\neq0,
		\]
		then the Hamiltonian system associated with \eqref{eq:hamiltonian-P3} possesses an unbounded solution. More precisely,
		if $\bigl[P_3\bigr]'(\theta_*)<0,$ then there exists a forward solution whose radial variable satisfies	$r(t)\rightarrow+\infty$ as $t\rightarrow+\infty;$ and if $\bigl[P_3\bigr]'(\theta_*)>0,$ then there exists a backward solution satisfying $r(t)\rightarrow+\infty$ as $t\rightarrow-\infty.$ In particular, in either case the corresponding solution of the original equation is unbounded.
	\end{theorem}

	The proof of this theorem is parallel to that of Theorem \ref{thm:main}. Now we give a brief sketch of the proof. For convenience, write
	\[
	A(\theta)=[P_3](\theta),
	\qquad
	B(\theta)=[P_4](\theta),
	\]
	and introduce
	\[
	u=r^{-1/2}.
	\]
	Thus, $u=0$ corresponds to $r=+\infty$. The Hamiltonian equations give
	\[
	\frac{du}{dt}
	=
	\frac12 A'(\theta)u^5
	+\frac12 B'(\theta)u^6
	+O(u^7)
	\]
	and
	\[
	\frac{d\theta}{dt}
	=
	-A(\theta)u^4
	-\frac32 B(\theta)u^5
	+O(u^6).
	\]
	In comparison with the case in which $[P_2]$ is the first nonvanishing
	term, the common degeneracy factor is now $u^4$. We therefore introduce
	the new time variable $s$ by
	\[
	\frac{ds}{dt}=u^4.
	\]
	For $u>0$, this is an orientation-preserving time transformation, and
	the desingularized system becomes
	\begin{equation}
		\begin{cases}
			\displaystyle
			\frac{du}{ds}
			=
			\frac12 A'(\theta)u
			+\frac12 B'(\theta)u^2
			+O(u^3),
			\\[2mm]
			\displaystyle
			\frac{d\theta}{ds}
			=
			-A(\theta)
			-\frac32 B(\theta)u
			+O(u^2).
		\end{cases}
		\label{eq:P3-desingularized-system}
	\end{equation}
	
	Let $\theta_*$ be a nondegenerate zero of $A$, namely,
	\[
	A(\theta_*)=0,
	\qquad
	\alpha:=A'(\theta_*)\neq0,
	\]
	and put
	\[
	\xi=\theta-\theta_*,
	\qquad
	B_*=B(\theta_*).
	\]
	The linearization of \eqref{eq:P3-desingularized-system} at
	$(u,\xi)=(0,0)$ is
	\[
	\frac{d}{ds}
	\begin{pmatrix}
		u\\ \xi
	\end{pmatrix}
	=
	\begin{pmatrix}
		\dfrac{\alpha}{2}&0\\[1mm]
		-\dfrac{3}{2}B_*&-\alpha
	\end{pmatrix}
	\begin{pmatrix}
		u\\ \xi
	\end{pmatrix}
	+
	O\!\left((|u|+|\xi|)^2\right).
	\]
	Its eigenvalues are
	\[
	\lambda_1=\frac{\alpha}{2},
	\qquad
	\lambda_2=-\alpha.
	\]
	Since $\alpha\neq0$, these eigenvalues have opposite signs. Hence
	$(u,\xi)=(0,0)$ is a hyperbolic saddle point.
	
	Assume first that
	\[
	\alpha=[P_3]'(\theta_*)<0.
	\]
	Then the eigendirection corresponding to $\lambda_1=\alpha/2<0$ has a
	nonzero $u$-component, and one of its two branches enters the physical region $u>0$.
	After introducing the linearizing coordinate
	\[
	z=\xi+\frac{B_*}{\alpha}u,
	\]
	the system has the standard saddle-point form
	\[
	\begin{cases}
		\displaystyle
		\frac{du}{ds}=\frac{\alpha}{2}u+F(u,z,t(s)),\\[2mm]
		\displaystyle
		\frac{dz}{ds}=-\alpha z+G(u,z,t(s)),
	\end{cases}
	\qquad
	F(u,z,t(s)),\,G(u,z,t(s))
	=
	O\!\left((|u|+|z|)^2\right).
	\]
	Although the linear part of the desingularized system has a
	hyperbolic saddle structure, the remainder still depends periodically
	on $t$, while $dt/ds=u^{-4}$ is singular at $u=0$. We therefore pass
	to the Poincar\'{e} map and repeat the discrete Lyapunov--Perron
	construction used in the nondegenerate $[P_2]$ case, with the
	corresponding higher-order polynomial weights. This yields a
	one-dimensional stable invariant curve in the physical region $u>0$
	of the form
	\[
	z=h(u),\qquad h(u)=O(u^2).
	\]
	In the original variables, this invariant curve satisfies
	\[
	\theta-\theta_*
	=
	-\frac{B(\theta_*)}{A'(\theta_*)}u
	+O(u^2)
	=
	-\frac{[P_4](\theta_*)}{[P_3]'(\theta_*)}u
	+O(u^2).
	\]
	
	A solution on the branch of the stable invariant curve contained in $u>0$
	satisfies
	\[
	u(s)\longrightarrow0
	\quad\text{as }\quad s\longrightarrow+\infty.
	\]
	It follows that
	\[
	r(s)=u(s)^{-2}\longrightarrow+\infty.
	\]
	Moreover,
	\[
	\frac{dt}{ds}=u(s)^{-4},
	\]
	and hence
	\[
	t(s)
	=
	t(0)+\int_0^s u(\sigma)^{-4}\,d\sigma
	\longrightarrow+\infty.
	\]
	Thus, the convergence $u(s)\to0$ corresponds to an orbit satisfying
	\[
	r(t)\longrightarrow+\infty
	\quad\text{as }\quad t\longrightarrow+\infty,
	\]
	rather than to a finite-time singularity.
	
	If the chosen nondegenerate zero satisfies
	\[
	[P_3]'(\theta_*)>0,
	\]
	the same construction applied in backward time produces an orbit
	unbounded in negative time.
	
	The proof is therefore analogous to the case in which
	$[P_2]$ has a nondegenerate zero. The only essential changes are that
	the leading Hamiltonian term is shifted from $r^{-1/2}$ to $r^{-1}$
	and that the desingularizing time transformation is correspondingly
	changed from
	\[
	\frac{ds}{dt}=u^3
	\qquad\text{to}\qquad
	\frac{ds}{dt}=u^4.
	\]
	After this modification, the desingularized system again possesses a
	saddle point, and the desired unbounded orbit is obtained from the
	corresponding stable invariant curve by an analogous discrete Lyapunov--Perron contraction argument.

	\section{Application}
	In this section, we investigate the boundedness and unboundedness of solutions to the equation
	\begin{equation}\label{eq:main}
		x''+4x^+-x^- = \cos kt,
		\qquad k\in\mathbb Z_+.
	\end{equation}

	We first prove a lemma showing that, in the high-energy region, the corresponding time-\(6\pi\) map is the identity for odd \(k\ge3\), whereas the corresponding time-\(3\pi\) map is the identity for \(k\equiv2\pmod 4\) with \(k\ge6\). Consequently, in these cases every sufficiently high-energy solution is periodic and, in particular, bounded.
	
	\begin{lemma}\label{periodic}
		Assume that $k\ge3, 4\nmid k,$ then, on the high-energy region, the successor dynamics is periodic.
		More precisely, if $k$ is odd, then $P_{6\pi}=\Id;$ and if $k\equiv2\pmod4$, then $P_{3\pi}=\Id.$ In either case,
		$P_{2\pi}^3=\Id$ for the usual time-$2\pi$ Poincar\'{e} map on the same high-energy region.
	\end{lemma}
	
	\noindent\textit{Proof.}	
	\textbf{Case I: odd $k\ge3$.} On $x<0$, equation \eqref{eq:main} becomes
	\[
	x''+x=\cos kt.
	\]
	Since $k\ne1$, a particular solution is
	\[
	p_-(t)=a_-\cos kt,
	\qquad a_-:=\frac{1}{1-k^2}.
	\]
	Let a negative-side segment start at a downward zero crossing
	\[
	x(\sigma)=0,
	\qquad x'(\sigma)=q<0.
	\]
	Then, for $s=t-\sigma$,
	\begin{equation}
		x(\sigma+s)
		=a_-\cos k(\sigma+s)
		-a_-\cos k\sigma\cos s
		+\bigl(q+a_-k\sin k\sigma\bigr)\sin s.
		\label{eq:negative-solution-odd}
	\end{equation}
	Since $k$ is odd,
	\[
	\cos k(\sigma+\pi)=-\cos k\sigma,
	\qquad
	\sin k(\sigma+\pi)=-\sin k\sigma.
	\]
	Putting $s=\pi$ in \eqref{eq:negative-solution-odd} gives
	\[
	x(\sigma+\pi)=0.
	\]
	Differentiating \eqref{eq:negative-solution-odd} gives
	\[
	x'(\sigma+s)
	=-a_-k\sin k(\sigma+s)
	+a_-\cos k\sigma\sin s
	+\bigl(q+a_-k\sin k\sigma\bigr)\cos s.
	\]
	Thus
	\[
	x'(\sigma+\pi)=-q.
	\]
	Therefore the negative half-axis produces the exact reflection
	\begin{equation}
		(\sigma,q)\longmapsto (\sigma+\pi,-q).
		\label{eq:negative-reflection-odd}
	\end{equation}
	For $|q|$ sufficiently large, the same explicit formula also shows that there
	is no other zero in $0<s<\pi$.  Indeed, the term $q\sin s$ dominates the bounded
	remaining terms on compact subintervals of $(0,\pi)$, and near $s=0$ and
	$s=\pi$ the signs are controlled by the nonzero derivatives $q$ and $-q$.
	
	On $x>0$,
	\[
	x''+4x=\cos kt.
	\]
	Since $k\ne2$, a particular solution is
	\[
	p_+(t)=a_+\cos kt,
	\qquad a_+:=\frac{1}{4-k^2}.
	\]
	Hence each positive-side segment can be written as
	\begin{equation*}
		x(t)=X_Z(t):=\operatorname{Re}(Ze^{2it})+a_+\cos kt,
	\end{equation*}
	where $Z\in\mathbb C$ is constant along that positive-side segment.  Equivalently,
	\begin{equation}
		Z=e^{-2it}\left[x(t)-a_+\cos kt
		-\frac{i}{2}\bigl(x'(t)+a_+k\sin kt\bigr)\right].
		\label{eq:Z-amplitude}
	\end{equation}
	
	Suppose a positive-side segment ends at a downward crossing
	\[
	x(\sigma)=0,
	\qquad x'(\sigma)=q<0.
	\]
	From \eqref{eq:Z-amplitude}, the old positive-side amplitude is
	\begin{equation*}
		Z=e^{-2i\sigma}\left[-a_+\cos k\sigma
		-\frac{i}{2}\bigl(q+a_+k\sin k\sigma\bigr)\right].
	\end{equation*}
	After the negative-side reflection \eqref{eq:negative-reflection-odd}, the next
	positive-side segment starts at
	\[
	x(\sigma+\pi)=0,
	\qquad x'(\sigma+\pi)=-q.
	\]
	Its new amplitude is therefore
	\begin{align*}
		Z_{\rm new}
		&=e^{-2i(\sigma+\pi)}\left[-a_+\cos k(\sigma+\pi)
		-\frac{i}{2}\bigl(-q+a_+k\sin k(\sigma+\pi)\bigr)\right] \\
		&=e^{-2i\sigma}\left[a_+\cos k\sigma
		+\frac{i}{2}\bigl(q+a_+k\sin k\sigma\bigr)\right]
		=-Z.
	\end{align*}
	Thus
	\begin{equation*}
		Z\longmapsto -Z.
	\end{equation*}
	Moreover,
	\begin{equation}
		X_{-Z}(t+\pi)=-X_Z(t).
		\label{eq:X-antiperiodic}
	\end{equation}
	Indeed,
	\[
	\begin{aligned}
		X_{-Z}(t+\pi)
		&=\operatorname{Re}\bigl(-Ze^{2i(t+\pi)}\bigr)
		+a_+\cos k(t+\pi)  \\
		&=-\operatorname{Re}(Ze^{2it})-a_+\cos kt
		=-X_Z(t),
	\end{aligned}
	\]
	where the oddness of $k$ has been used.
	
\begin{figure}[htb]
		\centering
		\includegraphics[width=0.8\linewidth]{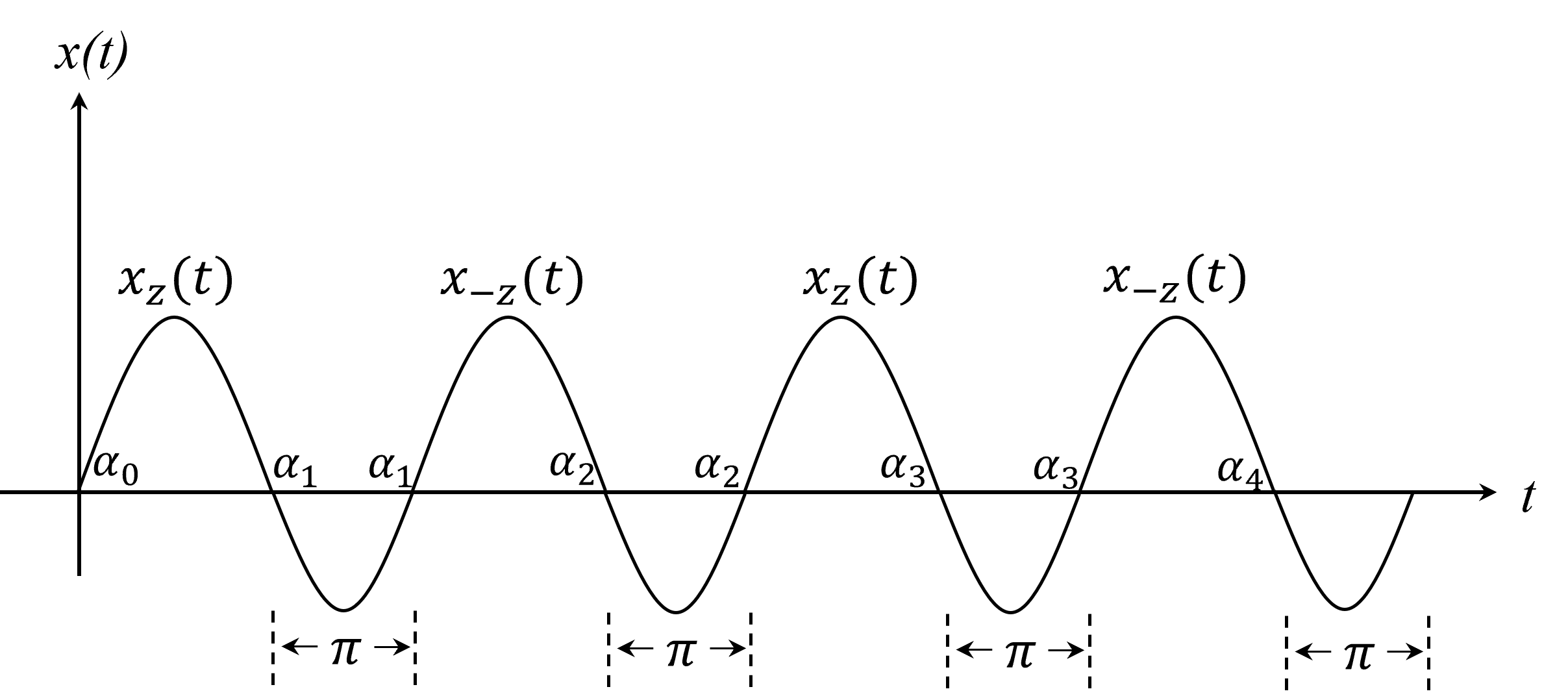}
		\caption{A schematic representation of \(x(t)\) for odd \(k\).}
		\label{phase}
	\end{figure}

	Let $|Z|$ be sufficiently large.  By the stability of simple zeros under small $C^1$ perturbations, $X_Z$ has
	exactly four simple zeros in one $2\pi$ period.  Write them as
	\[
	\alpha_0<\alpha_1<\alpha_2<\alpha_3<\alpha_4=\alpha_0+2\pi.
	\]
	Assume, without loss of generality, that
	\[
	X_Z(t)>0 \quad\text{for } t\in(\alpha_0,\alpha_1).
	\]
	Then the first positive-side segment has length $\alpha_1-\alpha_0$.
	After the following negative-side reflection, the positive amplitude becomes
	$-Z$.  By \eqref{eq:X-antiperiodic}, the next positive-side segment corresponds
	to the next sign interval of $X_Z$, and its length is
	$\alpha_2-\alpha_1$.  Repeating this argument, the four positive-side pieces have
	total length
	\[
	\begin{aligned}
		T_+
		&=(\alpha_1-\alpha_0)+(\alpha_2-\alpha_1)
		+(\alpha_3-\alpha_2)+(\alpha_4-\alpha_3) \\
		&=\alpha_4-\alpha_0=2\pi.
	\end{aligned}
	\]
	Each of the four negative-side pieces has exact length $\pi$, hence
	\[
	T_-=4\pi.
	\]
	The whole cycle therefore has length
	\[
	T=T_++T_-=6\pi.
	\]
	After the four negative-side reflections, the positive amplitude changes as
	\[
	Z\mapsto -Z\mapsto Z\mapsto -Z\mapsto Z.
	\]
	Thus the amplitude returns to its initial value.  See Figure~\ref{phase} for a schematic representation of the construction of the solution $x(t)$. In this figure, we glue together the solutions on the positive and negative half-axes.
	
	Since
	\[
	e^{2i(t+6\pi)}=e^{2it},
	\qquad
	\cos k(t+6\pi)=\cos kt,
	\]
	the solution and its velocity return to their initial values.  Hence
	\begin{equation*}
		P_{6\pi}=\Id
	\end{equation*}
	on the high-energy non-grazing region.
	
	\medskip
	\noindent\textbf{Case II: $k\equiv2\pmod4$, $k\ne2$.} Write \(k=4r+2\) with \(r\in\mathbb N\).
	
	On $x>0$,
	\[
	x''+4x=\cos kt,
	\]
	and, as before,
	\[
	p_+(t)=a_+\cos kt,
	\qquad a_+:=\frac{1}{4-k^2}.
	\]
	
	Let a positive-side segment start at an upward crossing
	\[
	x(\tau)=0,
	\qquad x'(\tau)=v>0.
	\]
	Then, for $s=t-\tau$,
	\begin{equation}
		x(\tau+s)
		=a_+\cos k(\tau+s)
		-a_+\cos k\tau\cos 2s
		+\frac12\bigl(v+a_+k\sin k\tau\bigr)\sin 2s.
		\label{eq:positive-solution-even2}
	\end{equation}
	Since $k\equiv2\pmod4$,
	\[
	\cos k\left(\tau+\frac{\pi}{2}\right)=-\cos k\tau,
	\qquad
	\sin k\left(\tau+\frac{\pi}{2}\right)=-\sin k\tau.
	\]
	Putting $s=\pi/2$ in \eqref{eq:positive-solution-even2} gives
	\[
	x\left(\tau+\frac{\pi}{2}\right)=0.
	\]
	Differentiating \eqref{eq:positive-solution-even2} gives
	\[
	x'(\tau+s)
	=-a_+k\sin k(\tau+s)
	+2a_+\cos k\tau\sin2s
	+\bigl(v+a_+k\sin k\tau\bigr)\cos2s.
	\]
	Thus
	\[
	x'\left(\tau+\frac{\pi}{2}\right)=-v.
	\]
	Therefore the positive half-axis gives the exact reflection
	\begin{equation}
		(\tau,v)\longmapsto\left(\tau+\frac{\pi}{2},-v\right).
		\label{eq:positive-reflection-even2}
	\end{equation}
	For $v$ sufficiently large, the explicit formula also implies that no extra
	zero appears in $0<s<\pi/2$.
	
	On $x<0$,
	\[
	x''+x=\cos kt.
	\]
	Since $k\ne1$, a particular solution is
	\[
	p_-(t)=a_-\cos kt,
	\qquad a_-:=\frac{1}{1-k^2}.
	\]
	Hence each negative-side segment can be written as
	\begin{equation*}
		x(t)=Y_W(t):=\operatorname{Re}(We^{it})+a_-\cos kt,
	\end{equation*}
	where $W\in\mathbb C$ is constant along that negative-side segment.  Equivalently,
	\begin{equation*}
		W=e^{-it}\left[x(t)-a_-\cos kt
		-i\bigl(x'(t)+a_-k\sin kt\bigr)\right].
	\end{equation*}
	
	Suppose a negative-side segment ends at an upward crossing
	\[
	x(\tau)=0,
	\qquad x'(\tau)=v>0.
	\]
	Then
	\begin{equation*}
		W=e^{-i\tau}\left[-a_-\cos k\tau
		-i\bigl(v+a_-k\sin k\tau\bigr)\right].
	\end{equation*}
	After the positive-side reflection \eqref{eq:positive-reflection-even2}, the
	next negative-side segment starts at
	\[
	x\left(\tau+\frac{\pi}{2}\right)=0,
	\qquad
	x'\left(\tau+\frac{\pi}{2}\right)=-v.
	\]
	The new negative-side amplitude is
	\begin{align*}
		W_{\rm new}
		&=e^{-i(\tau+\pi/2)}
		\left[-a_-\cos k\left(\tau+\frac{\pi}{2}\right)
		-i\left(-v+a_-k\sin k\left(\tau+\frac{\pi}{2}\right)\right)\right] \\
		&=-ie^{-i\tau}
		\left[a_-\cos k\tau+i\bigl(v+a_-k\sin k\tau\bigr)\right]
		=iW.
	\end{align*}
	Thus
	\begin{equation*}
		W\longmapsto iW.
	\end{equation*}
	Moreover,
	\begin{equation}
		Y_{iW}\left(t+\frac{\pi}{2}\right)=-Y_W(t).
		\label{eq:Y-antiperiodic}
	\end{equation}
	Indeed,
	\[
	\begin{aligned}
		Y_{iW}\left(t+\frac{\pi}{2}\right)
		&=\operatorname{Re}\left(iW e^{i(t+\pi/2)}\right)
		+a_-\cos k\left(t+\frac{\pi}{2}\right) \\
		&=\operatorname{Re}(-We^{it})-a_-\cos kt
		=-Y_W(t).
	\end{aligned}
	\]
	
\begin{figure}[htb]
		\centering
		\includegraphics[width=0.7\linewidth]{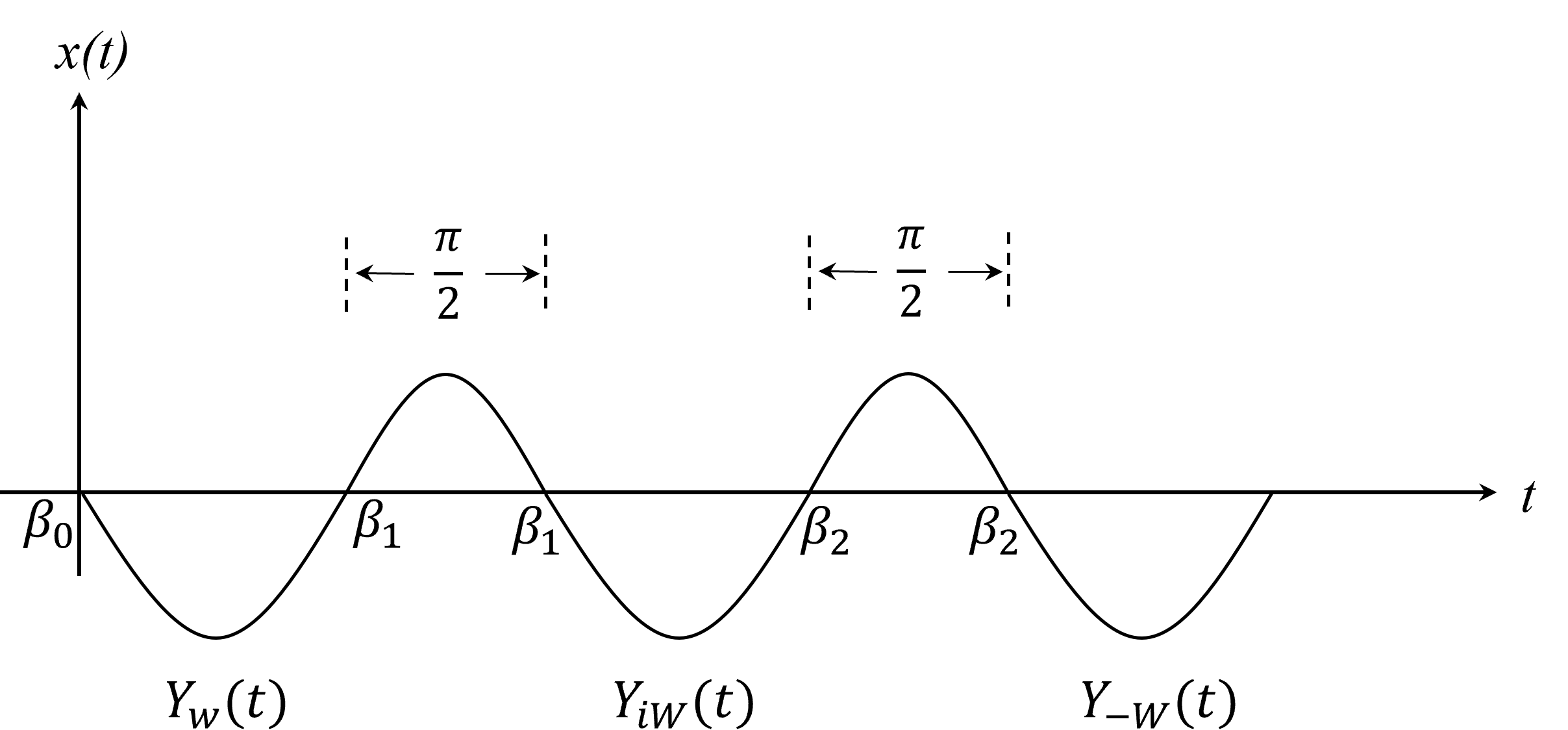}
		\caption{A schematic representation of \(x(t)\) for $k\equiv2\pmod4$, $k\ne2$.}
		\label{phase2}
	\end{figure}

	Let $|W|$ be sufficiently large.  By the stability of simple zeros under small $C^1$ perturbations, $Y_W$ has
	exactly two simple zeros in one $2\pi$ period.  Denote them by
	\[
	\beta_0<\beta_1<\beta_2=\beta_0+2\pi.
	\]
	Assume, without loss of generality, that
	\[
	Y_W(t)<0 \quad\text{for } t\in(\beta_0,\beta_1).
	\]
	Then the first negative-side segment has length $\beta_1-\beta_0$.
	After the following positive-side reflection, the negative-side amplitude
	becomes $iW$.  By \eqref{eq:Y-antiperiodic}, the second negative-side segment
	corresponds to the next sign interval of $Y_W$, and has length
	$\beta_2-\beta_1$.  Therefore the total negative-side time is
	\[
	T_-=(\beta_1-\beta_0)+(\beta_2-\beta_1)=2\pi.
	\]
	There are two positive-side reflections, each of length $\pi/2$, hence
	\[
	T_+=\pi.
	\]
	Thus the full cycle has length
	\[
	T=T_-+T_+=3\pi.
	\]
	After two positive-side reflections, the negative-side amplitude changes as
	\[
	W\mapsto iW\mapsto -W.
	\]
	The time has increased by $3\pi$. See Figure~\ref{phase2} for a schematic representation of the construction of the solution $x(t)$. In this figure, we glue together the solutions on the positive and negative half-axes.

	Since
	\[
	e^{i(t+3\pi)}=-e^{it},
	\qquad
	\cos k(t+3\pi)=\cos kt
	\quad (k\equiv2\pmod4),
	\]
	we get
	\[
	Y_{-W}(t+3\pi)=Y_W(t),
	\qquad
	Y_{-W}'(t+3\pi)=Y_W'(t).
	\]
	Therefore the solution and velocity return to their initial values, and hence
	\begin{equation*}
		P_{3\pi}=\Id
	\end{equation*}
	on the high-energy region.
	
	\hfill\qedsymbol

	In fact, the preceding derivation can be generalized to the equation
	$x''+p^2x^{+}-q^2x^{-}=\cos jt$, where $p,q\in\mathbb Z_{+},j\in\mathbb Z_{+}$. Assume that $j\neq p$, $j\neq q$. If either $j=p(2N+1)$ for some $N\in\mathbb Z$ or $j=q(2N+1)$ for some $N\in\mathbb Z$, then the corresponding high-energy successor map is periodic on the set of nondegenerate crossing orbits. More precisely, there exists \(R>0\) such that the following conclusions hold for every solution whose initial energy is greater than \(R\) and whose crossings of \(x=0\) are all transverse:
	
	\begin{enumerate}
		\item If
		\[
		j=q(2N+1),\qquad j\neq q,
		\]
		then the negative half-line equation possesses an exact half-period reflection. Consequently, the time-\(T_{q}\) map is the identity on the high-energy nondegenerate crossing region, where
		\[
		T_{q}
		=
		2\pi+\frac{2p\pi}{q}
		=
		\frac{2\pi(p+q)}{q}.
		\]
		
		\item If
		\[
		j=p(2N+1),\qquad j\neq p,
		\]
		then the positive half-line equation possesses an exact half-period reflection. Consequently, the time-\(T_{p}\) map is the identity on the high-energy nondegenerate crossing region, where
		\[
		T_{p}
		=
		2\pi+\frac{2q\pi}{p}
		=
		\frac{2\pi(p+q)}{p}.
		\]
	\end{enumerate}
	
	We now apply the unbounded-motion criteria established in the preceding sections to the asymmetric oscillator \eqref{eq:main}.
	
	\medskip
	\noindent\textbf{Proof of Theorem~\ref{th2}.} For $a=4$ and $b=1$, the parameter $\omega_0$ introduced in the normal-form construction is
	\[
	\omega_{0}
	=\frac{1}{2}\left(\frac{1}{\sqrt a}+\frac{1}{\sqrt b}\right)
	=\frac34.
	\]
	Consequently, in the notation of the normal-form construction,
	\[
	n=3,
	\qquad
	m=4,
	\qquad
	\rho=\sqrt{\frac{2m}{an}}=\frac{\sqrt6}{3}.
	\]
	The corresponding constants are
	\[
	\kappa
	=\rho m\left(\frac nm\right)^{3/2}
	=\frac{3\sqrt2}{2},
	\qquad
	\beta=\frac{n^{2}\rho^{2}}{2m}=\frac34,
	\qquad
	\delta
	=\frac{m\rho^{3}}{8}\left(\frac nm\right)^{5/2}
	=\frac{3\sqrt2}{32}.
	\]
	
	Let
	\[
	c(t)=C(3t),
	\qquad
	F_k(\theta,t)=\cos k(\theta+3t),
	\]
	where $C$ is the asymmetric cosine associated with $a=4$ and $b=1$.
	The function $c$ is $\pi/2$-periodic.
	
	For $k=2$, direct substitution into the first averaged coefficient and
	piecewise integration over the smoothness intervals of $C(3t)$ give
	\begin{equation*}
		[P_1](\theta)
		=-\frac{61}{192}+\frac{1}{32}\cos4\theta.
	\end{equation*}
	Thus $[P_1]$ is not constant. The first-order unbounded-motion criterion
	therefore applies and yields that the equation
	\[
	\ddot{x}+4x^{+}-x^{-}=\cos2t
	\]
	possesses an unbounded solution.

	In the case $k=1$ the first two relevant averaged terms degenerate, and it is
	necessary to proceed to the coefficient $[P_3]$. Indeed, since $c(t)$ is $\pi/2$-periodic, its Fourier expansion contains only temporal
	frequencies $4j$, $j\in\mathbb Z$. On the other hand,
	\[
	F_1(\theta,t)
	=\frac12e^{i(\theta+3t)}
	+\frac12e^{-i(\theta+3t)}.
	\]
	Thus the two Fourier modes of $F_1$ have bidegrees
	\[
	(p,q)=(1,3)
	\quad\text{and}\quad
	(p,q)=(-1,-3),
	\]
	where $p$ denotes the $\theta$-frequency and $q$ the $t$-frequency. Since $Q$ is chosen as the zero-mean periodic primitive of $c(t)F_1(\theta,t)$, it contains no Fourier mode with zero $t$-frequency. Moreover, its nonzero Fourier modes satisfy the same selection rule as those of $c(t)F_1(\theta,t)$. Consequently, every Fourier monomial produced in the normal-form calculation satisfies the
	selection rule
	\begin{equation}\label{eq:selection}
		q-3p\in4\mathbb Z.
	\end{equation}
	
	The function $P_1$ is quadratic in the first-order quantities. Hence its
	possible angular frequencies are
	\[
	p\in\{-2,0,2\}.
	\]
	After averaging with respect to $t$, only terms with $q=0$ remain. By
	\eqref{eq:selection}, this requires $3p\in4\mathbb Z$, and therefore
	$p\in4\mathbb Z$. The only admissible value is $p=0$. It follows that
	\begin{equation*}
		[P_1](\theta)\equiv c_1,
	\end{equation*}
	where $c_1$ is independent of $\theta$.
	
	Now we compute $[P_2](\theta)$. Let
	\[
	Q(\theta,t)=\int_{0}^{t}c(s)F_1(\theta,s)\,ds
	-\frac{1}{2\pi}\int_{0}^{2\pi}\left(\int_{0}^{\tau}c(s)F_1(\theta,s)\,ds\right)d\tau,
	\qquad
	Q_j=\partial_\theta^{j}Q.
	\]
	The second coefficient has the form
	\begin{equation*}
		P_2
		=G_1+\frac{\kappa Q}{2}\,\partial_\theta G_0,
	\end{equation*}
	where
	\begin{align*}
		G_0
		&=\beta c^2F_1^2
		-\frac{\kappa^2}{2}cF_1Q_1,\\
		G_1
		&=\delta c^3F_1^3
		-\frac{\kappa^3}{8}cF_1Q_1^2.
	\end{align*}
	Every term in $P_2$ is cubic in the first-order oscillatory quantities. Its
	possible angular frequencies are consequently
	\[
	p\in\{-3,-1,1,3\}.
	\]
	If such a term survives time averaging, then $q=0$. The selection rule
	\eqref{eq:selection} again implies $p\in4\mathbb Z$. Since
	\[
	\{-3,-1,1,3\}\cap4\mathbb Z=\varnothing,
	\]
	none of these modes can survive. Hence
	\begin{equation*}
		[P_2](\theta)\equiv0.
	\end{equation*}
	
	Then we consider $[P_3]$. The third coefficient is
	\begin{equation}\label{eq:P3-form}
		\begin{aligned}
			P_3={}&G_2
			+\frac{\kappa Q}{2}\,\partial_\theta G_1
			+\frac{\kappa^2QQ_1}{4}\,\partial_\theta G_0
			+\frac{\kappa^2Q^2}{8}\,\partial_\theta^2G_0,
		\end{aligned}
	\end{equation}
	where
	\[
	G_2
	=-\frac{\kappa^4}{16}cF_1Q_1^3
	+\frac{\kappa\delta}{2}c^3F_1^3Q_1.
	\]
	The fourth-order frequency selection rule shows that
	\begin{equation}\label{eq:P3-fourier}
		[P_3](\theta)
		=A_0+A_c\cos4\theta+A_s\sin4\theta.
	\end{equation}
	We next show analytically that $A_0=A_s=0$.
	
	Set
	\[
	\mathcal A(\theta,t)=\frac{\kappa}{2}Q(\theta,t).
	\]
	Since $\partial_\theta\mathcal A=\kappa Q_1/2$, one has
	\[
	(\partial_\theta\mathcal A)G_1=G_2.
	\]
	Moreover,
	\[
	\mathcal A\,\partial_\theta\mathcal A
	=\frac{\kappa^2QQ_1}{4},
	\qquad
	\frac{\mathcal A^2}{2}=\frac{\kappa^2Q^2}{8}.
	\]
	Consequently, \eqref{eq:P3-form} can be rewritten as the exact angular
	derivative
	\begin{equation}\label{eq:P3-total-derivative}
		P_3(\theta,t)
		=\partial_\theta\left(
		\frac{\kappa Q}{2}G_1
		+\frac{\kappa^2Q^2}{8}\,\partial_\theta G_0
		\right).
	\end{equation}
	It follows that
	\[
	\int_{0}^{2\pi}P_3(\theta,t)\,d\theta=0
	\]
	for every $t$. Averaging in $t$ and comparing with
	\eqref{eq:P3-fourier} gives
	\begin{equation}\label{eq:A0-zero}
		A_0=0.
	\end{equation}
	
	To eliminate the sine term, observe that $c$ is even and
	\[
	F_1(-\theta,-t)=F_1(\theta,t).
	\]
	The definition of $Q$ then yields
	\[
	Q(-\theta,-t)=-Q(\theta,t),
	\qquad
	Q_1(-\theta,-t)=Q_1(\theta,t).
	\]
	Using \eqref{eq:P3-total-derivative}, one obtains
	\[
	P_3(-\theta,-t)=P_3(\theta,t).
	\]
	Time averaging over a complete period therefore gives
	\[
	[P_3](-\theta)=[P_3](\theta).
	\]
	Thus $[P_3]$ is even, and comparison with
	\eqref{eq:P3-fourier} gives
	\begin{equation}\label{eq:As-zero}
		A_s=0.
	\end{equation}
	Combining \eqref{eq:P3-fourier}, \eqref{eq:A0-zero}, and
	\eqref{eq:As-zero}, we obtain the exact identity
	\begin{equation*}
		[P_3](\theta)=A_c\cos4\theta.
	\end{equation*}

	Numerical piecewise integration gives
	\begin{equation*}
		A_c=[P_3](0)
		\approx-0.08875874269.
	\end{equation*}
	Since $A_c\neq0$, $[P_3](\theta)=A_c\cos4\theta$ has nondegenerate zeros. Therefore, by Theorem~\ref{thm:P3-nondegenerate-zero}, equation \eqref{eq:main} with $k=1$ possesses an unbounded solution.
	
	Finally, suppose that \(k\ge3\) and \(4\nmid k\). By Lemma ~\ref{periodic}, the successor dynamics is periodic in the high-energy region. Hence an orbit entering this region cannot escape to infinity. Therefore, every solution of equation \eqref{eq:main} is bounded.
	
\hfill\qedsymbol


\section*{Data availability}
No data was used for the research described in the article.

\section*{References}
	{\small
		\rm
		
	}
	
\end{document}